\documentclass[a4paper,12pt]{amsart}

\usepackage[T1]{fontenc}
\usepackage[utf8]{inputenc}
\usepackage[english]{babel}
\usepackage{fullpage}
\usepackage[colorlinks]{hyperref}
\usepackage{amssymb,amsfonts,amsxtra,amsmath}
\usepackage{stmaryrd}
\usepackage[all,cmtip]{xy}
\usepackage{paralist}
\usepackage[shortlabels]{enumitem}
\usepackage{tikz,tikz-3dplot}
\usepackage[small]{caption}

\theoremstyle{plain}
\newtheorem{thm}{Theorem}[subsection]
\newtheorem{cor}[thm]{Corollary} 
\newtheorem{prp}[thm]{Proposition} 
\newtheorem{lem}[thm]{Lemma} 
\theoremstyle{definition}
\newtheorem{dfn}[thm]{Definition}

\theoremstyle{remark}
\newtheorem{rmk}[thm]{Remark}

\newtheorem{exa}[thm]{Example}

\newcommand{\NN}{\mathbb{N}}
\newcommand{\ZZ}{\mathbb{Z}}
\newcommand{\ZI}{\mathbb{Z}_\infty}
\newcommand{\CC}{\mathbb{C}}

\newcommand{\C}{\mathcal{C}}
\newcommand{\E}{\mathcal{E}}
\newcommand{\F}{\mathcal{F}}
\newcommand{\G}{\mathcal{G}}
\newcommand{\K}{\mathcal{K}}

\newcommand{\Q}{\mathcal{Q}}
\newcommand{\V}{\mathfrak{V}}
\newcommand{\RFI}{\mathfrak{R}}
\newcommand{\GSI}{\mathfrak{G}}

\newcommand{\ee}{\mathbf{e}}

\newcommand{\mm}{\mathfrak{m}}
\newcommand{\nn}{\mathfrak{n}}
\newcommand{\pp}{\mathfrak{p}}

\newcommand{\rr}{\mathfrak{r}}

\newcommand{\reg}{\mathrm{reg}}

\newcommand{\ol}{\overline}
\newcommand{\wh}{\widehat}
\newcommand{\wt}{\widetilde}

\newcommand{\ideal}[1]{{\left\langle#1\right\rangle}}

\newcommand{\onto}{\twoheadrightarrow}

\DeclareMathOperator{\depth}{depth}

\DeclareMathOperator{\Hom}{Hom}
\DeclareMathOperator{\Ext}{Ext}
\DeclareMathOperator{\Spec}{Spec}

\DeclareMathOperator{\Max}{Max}

\numberwithin{equation}{section}

\begin{document}

\title{Duality on value semigroups}

\author[Ph.~Korell]{Philipp Korell}
\address{Ph.~Korell\\
Department of Mathematics\\
TU Kaiserslautern\\
67663 Kaiserslautern\\
Germany}
\email{\href{mailto:korell@mathematik.uni-kl.de}{korell@mathematik.uni-kl.de}}

\author[M.~Schulze]{Mathias Schulze}
\address{M.~Schulze\\
Department of Mathematics\\
TU Kaiserslautern\\
67663 Kaiserslautern\\
Germany}
\email{\href{mailto:mschulze@mathematik.uni-kl.de}{mschulze@mathematik.uni-kl.de}}

\author[L.~Tozzo]{Laura Tozzo}
\address{L.~Tozzo\\
Department of Mathematics\\
TU Kaiserslautern\\
67663 Kaiserslautern\\
Germany}
\email{\href{mailto:tozzo@mathematik.uni-kl.de}{tozzo@mathematik.uni-kl.de}}

\thanks{The research leading to these results has received funding from the People Programme (Marie Curie Actions) of the European Union's Seventh Framework Programme (FP7/2007-2013) under REA grant agreement n\textsuperscript{o} PCIG12-GA-2012-334355.}




\begin{abstract}
We establish a combinatorial counterpart of the Cohen--Macaulay duality on a class of curve singularities which includes algebroid curves.
For such singularities the value semigroup and the value semigroup ideals of all fractional ideals satisfy axioms that define so-called good semigroups and good semigroup ideals. 
We prove that each good semigroup admits a canonical good semigroup ideal which gives rise to a duality on good semigroup ideals.
We show that the Cohen--Macaulay duality and our good semigroup duality are compatible under taking values.
\end{abstract}

\maketitle
\tableofcontents

\section{Introduction}

Value semigroups of curve singularities have been studied intensively for decades.
Waldi~\cite{Wal72,Wal00} showed that any plane algebroid curve is determined by its value semigroup up to equivalence in the sense of Zariski.
The value semigroup thus determines the topological type for any plane complex curve singularity.

Kunz~\cite{Kun70} showed that an analytically irreducible and residually rational local ring $R$ is Gorenstein if and only if its (numerical) value semigroup $\Gamma_R$ is \emph{symmetric}. 
J\"ager~\cite{Jag77} used the symmetry condition to define a semigroup ideal $K^0$ such that (suitably normalized) canonical ideals $\K$ of $R$ are characterized by having value semigroup ideal $\Gamma_\K=K^0$.

Waldi~\cite{Wal72} was the first to describe a symmetry property of the value semigroup for plane algebroid curves with two branches.
Delgado~\cite{Del87,Del88} then made the step to general algebroid curves proving an analog of Kunz's result.
Later Campillo, Delgado, and Kiyek \cite{CDK94} relaxed the hypotheses to include analytically reduced and residually rational local rings $R$ with infinite residue field.

D'Anna~\cite{DAn97} extended J\"ager's approach under the preceding hypotheses.
He turned Delgado's symmetry definition into an explicit formula for a semigroup ideal $K^0$ (see Definition~\ref{46}) such that any (suitably normalized) fractional ideal $\K$ of $R$ is canonical if and only if $\Gamma_\K=K^0$.
In the process he studied axioms satisfied by value semigroup ideals which lead to the notion of a \emph{good semigroup ideal} (see Definition~\ref{26}).

Barucci, D'Anna, and Fr{\"o}berg~\cite{BDF00} studied some more special classes of rings such as almost Gorenstein rings, Arf rings, and rings of small multiplicity in relation with their value semigroups.
Their setup includes the case of semilocal rings.
Notably they found an example of a good semigroup which is not the value semigroup of any ring.

Recently Pol~\cite[Thm.~5.2.1]{Pol16} gave an explicit formula for the value semigroup ideal of the dual of a fractional ideal for Gorenstein algebroid curves.

\medskip

In this paper, we extend and unify D'Anna's and Pol's results for a general class of rings $R$ that we call \emph{admissible} (see Definition~\ref{88}).
We show that any good semigroup admits a \emph{canonical semigroup ideal} $K$ that is defined by a simple maximality property (see Definition~\ref{44}).
Equivalently, such a $K$ induces a duality $E\mapsto K-E$ on good semigroup ideals (see Theorem~\ref{45}).
This means that
\[
K-(K-E)=E
\]
for all good semigroup ideals.
It turns out that our canonical semigroup ideals are exactly the translations of D'Anna's $K^0$.
In particular, D'Anna's characterization of canonical ideals in terms of their value semigroup ideals persists for admissible rings (see Theorem~\ref{111}).
For any canonical ideal $\K$ of $R$ we show that 
\[
\Gamma_{\K:\E}=\Gamma_\K-\Gamma_\E
\]
for all regular fractional ideals $\E$ of $R$ (see Theorem~\ref{74}).
This means that there is a commutative diagram
\begin{equation*}
\xymatrixrowsep{0.7pc}
\xymatrixcolsep{1.5pc}
\xymatrix{
\left\{\vcenter{\txt{regular fractional \\ ideals of $R$}}\right\} \ar[rr]^-{\E\mapsto\K:\E}\ar[dd]_-{\E\mapsto\Gamma_\E} && \left\{\vcenter{\txt{regular fractional \\ ideals of $R$}}\right\}\ar[dd]^-{\E\mapsto\Gamma_\E}\\
& \circlearrowleft & \\
\left\{\vcenter{\txt{good semigroup \\ ideals of $\Gamma_R$}} \right\}\ar[rr]_-{E\mapsto\Gamma_\K-E} &&\left\{\vcenter{\txt{good semigroup \\ ideals of $\Gamma_R$}} \right\}
}
\end{equation*}
relating the Cohen--Macaulay duality $\E\mapsto\K:\E$ on $R$ to our good semigroup duality $E\mapsto K-E$ on $\Gamma_R$ for $K=\Gamma_\K$.

An important tool to prove the commutativity of the above diagram is the \emph{distance} $d(F\backslash E)$ between two good semigroup ideals $E\subset F$ (see Definition~\ref{113}).
It plays the role of the length $\ell_R(\F/\E)$ of the quotient of two fractional ideals $\E\subset\F$ on the semigroup side.
In fact, the two quantities agree in the case where $E=\Gamma_\E$ and $F=\Gamma_\F$ (see Proposition~\ref{39}), that is,
\[
\ell_R(\F/\E)=d(\Gamma_\F\backslash \Gamma_\E).
\]
D'Anna~\cite[2.7~Prop.]{DAn97} stated that $d(F\backslash E)=0$ is equivalent to $E=F$, which implies $\E=\F$ in the preceding case.
We give a proof of this crucial fact (see Proposition~\ref{37}).

Before approaching these main results, we review the definition of value semigroups and their ideals and give a detailed account of their compatibility with localization and completion (see \S\ref{77}).

\section{Preliminaries}\label{85}

All rings under consideration are commutative and unitary.
For a ring $R$ we denote by $\Max(R)$ the set of its maximal ideals.
We call a one-dimensional Noetherian ring $R$ Cohen--Macaulay if $\depth(R_\mm)=1$ for all $\mm\in\Max(R)$.

For an $R$-module $M$ we write $\ell_R(M)$ for its length and $\wh M$ for its completion at the Jacobson radical of $R$.
By $\ee_i$ we denote the $i$th unit vector of a free module.

The \emph{total ring of fractions} $Q_R$ of a ring $R$ is the localization of $R$ at the set $R^\reg$ of all regular elements of $R$.
More generally, we set $S^\reg:=S\cap Q_R^\reg$ for any subset $S\subset Q_R$.
Note that $R^\reg=R\cap Q_R^\reg$.
We denote by $\ol R$ the integral closure of $R$ in $Q_R$.
If $R$ is reduced, then it coincides with the \emph{normalization} of $R$.

We abbreviate $\F:\E:=\F:_{Q_R}\E$ for any subsets $\E,\F\subset Q_R$.
We collect some trivial properties of this colon operation for future reference.


\begin{rmk}\label{108}
Let $x\in Q_R^\reg$ and $\E,\E',\F,\F',\G$ be $R$-submodules of $Q_R$.
Then
\begin{enumerate}[label=(\alph*), ref=\alph*]
\item\label{108a} $(\G:\F):\E=\G:(\F\E)$,
\item\label{108b} $(x\E):\F=x(\E:\F)=\E:(x^{-1}\F)$,
\item\label{108c} $\E:\F'\subset\E:\F\subset\E':\F$ if $\E \subset \E'$ and $\F \subset \F'$, and
\item\label{108d} $\E:\F=(\E:R'):\F$ if $R\subset R'\subset Q_R$ is a ring extension and $\F$ an $R'$-module. 
\end{enumerate}
\end{rmk}

\subsection{Fractional ideals}\label{83}

Fractional ideals play a central role in our considerations.
Here we summarize the properties we shall use. 
Let $R$ be a ring.


\begin{dfn}\label{106}\
\begin{enumerate}[label=(\alph*), ref=\alph*]
\item\label{106a} An $R$-submodule $\E$ of $Q_R$ is called \emph{regular} if $\E^\reg\ne\emptyset$ or, equivalently, $Q_R\E=Q_R$.

\item\label{106b} An $R$-submodule $\E\subset Q_R$ such that $r\E\subset R$ for some $r\in R^\reg$ is called a \emph{fractional ideal} (of $R$).
We denote by $\RFI_R$ the set of regular fractional ideals of $R$.

\item\label{106e} An $R$-submodule $\E$ of $Q_R$ is \emph{invertible} if $\E\F=R$ for some $R$-submodule $\F$ of $Q_R$.
We denote by $\RFI_R^*$ the set of all invertible $R$-submodules of $Q_R$.

\item\label{106d} The \emph{conductor} of a fractional ideal $\E$ of $R$ is $\C_\E := \E : \ol R\subset\E$.
\end{enumerate}
\end{dfn}


\begin{rmk}\label{134}
The fractional ideals of a Noetherian ring $R$ are the finitely generated $R$-submodules of $Q_R$.
If $R$ is a one-dimensional Cohen--Macaulay ring, then any $\F\in\RFI_R$ is a faithful maximal Cohen-Macaulay module of $R$.
\end{rmk}


The set $\RFI_R$ is a (commutative) monoid under product of ideals.
If $\E\subset Q_R$ is an invertible $R$-submodule, then it is regular and finitely generated, and its inverse is uniquely determined as $\F=\E^{-1}=R:\E$ (see \cite[Ch.~II, (2.1) Rem.~(3) and (2.2) Prop.~(1), (2)]{KV04}).
In particular, the (abelian) group $\RFI_R^*$ of all invertible $R$-submodule of $Q_R$ is a submonoid of $\RFI_R$.
In the case where $R$ is semilocal, all elements of $\RFI_R^*$ are principal fractional ideals (see \cite[Ch.~II, (2.2) Prop.~(3)]{KV04}).


In the following we summarize the relation of the colon operation with the Hom functor and flat base change and well-known properties of completion.


\begin{lem}\label{2}
Let $R$ be a ring.

\begin{enumerate}[label=(\alph*), ref=\alph*]

\item\label{2a} For $\E,\F\in\RFI_R$, $\F:\E\in\RFI_R$, and there is a canonical isomorphism
\[
\F:\E\to\Hom_R(\E,\F),\quad x\mapsto(y\mapsto xy),
\]
of $R$-modules compatible with multiplication in $Q_R$ and composition of homomorphisms.
Iterating yields a commutative diagram of canonical maps
\[
\xymatrix{
\E \ar[r]\ar[dr] & \F:(\F:\E) \ar[d]^-\cong\\
& \Hom_R(\Hom_R(\E,\F),\F).
}
\]

\item\label{2b} Any flat ring homomorphism $\varphi\colon R\to R'$ induces a ring homomorphism
\[
\wt\varphi\colon Q_R\to Q_{R'}.
\]
If $\varphi$ is injective, then also $\wt\varphi$ is injective, and $\E R':=\wt\varphi(\E)R'\cong\E\otimes_RR'$ for any $R$-submodule $\E$ of $Q_R$.

\item\label{2c} If $R\to R'$ is flat and $\E,\F\in\RFI_R$, then $\E\otimes_RR'\cong\E R'\in\RFI_{R'}$, and there is a commutative diagram
\[
\xymatrix{
& Q_R\ar[d]\ar[dl]_-{\wt\varphi} & \ar@{_(->}[l] \F:\E\ar[d]\ar[r]_-\cong & \Hom_R(\E,\F)\ar[d]\\
Q_{R'} & Q_RR'\ar@{_(->}[l] & (\F:\E)R'\ar@{_(->}[l]\ar@{=}[d]\ar[r]_-\cong & \Hom_R(\E,\F)\otimes_RR'\ar[d]^-\cong\\
& & \ar@{_(->}[ull] \F R':\E R'\ar[r]_-\cong & \Hom_{R'}(\E R',\F R').
}
\]

\item\label{2d} If $\colon R\to R'$ is faithfully flat (and hence injective), then $\E R'\cap Q_R=\E$ and $(\E\cap\F)R'=\E R'\cap\F R'$ for any $R$-submodules $\E$ and $\F$ of $Q_R$.

\end{enumerate}
\end{lem}

\begin{proof}
See \cite[Lem.~2.1 and 2.3]{HK71} and \cite[Ch.~I, \S3, no.~5, Prop.~10]{Bou61}.
\end{proof}


\begin{lem}\label{5}
Let $R$ be a Noetherian ring.
\begin{enumerate}[label=(\alph*), ref=\alph*]
\item\label{5a} The ring extension $R\to\wh R$ is faithfully flat.
\item\label{5b} If $\E$ is finitely generated, then $\E\wh R=\wh\E$.
\item\label{5e} If $R$ is semilocal, then $\wh R = \prod_{\mm \in \Max(R)}\wh{R_\mm}$, where $\wh{R_\mm}=\wh R_\mm=\wh R_{\wh\mm}$ are local rings.
\item\label{5f} If $R$ is semilocal and $R\subset R'$ is a finite ring extension, then $R'\otimes_R\wh R=\wh{R'}$.
\end{enumerate}
\end{lem}

\begin{proof}
See \cite[Thms.~8.7, 8.14, 8.15]{Mat89} and \cite[(16.8) Thm.]{Nag62}.
To see $\wh R_\mm=\wh R_{\wh\mm}$ in \eqref{5e} note that $\mm\wh R=\wh\mm$ by \eqref{5b}, and hence $\mm=\wh\mm\cap R$ by \eqref{5a} and Lemma~\ref{2}.\eqref{2d}.
\end{proof}


In our main case of interest, regular fractional ideals are in bijection under completion.


\begin{lem}\label{6}
Let $R$ be a one-dimensional local Cohen--Macaulay ring.
Then $Q_R\wh R=Q_{\wh R}$, and there is an inclusion preserving group isomorphism
\begin{align*}
\RFI_R &\to \RFI_{\wh R},\\
\E &\mapsto\wh\E, \\
\quad\F \cap Q_R &\mapsfrom\F.
\end{align*}
\end{lem}

\begin{proof}
See \cite[Ch.~II, (2.4)]{KV04}, \cite[Lem.~2.11]{HK71} and Lemmas~\ref{2}.\eqref{2c} and \eqref{2d} and \ref{5}.\eqref{5a} and \eqref{5b}.
\end{proof}


The following result will serve to eliminate the ambiguity of canonical ideals.


\begin{lem}\label{7}
Let $R=(R,\mm)$ be a local Noetherian ring, $R'\subset Q_R$ a finite extension ring of $R$ with $\lvert R/\mm \rvert \ge \lvert \Max(R') \rvert$, and $\E\in\RFI_R$ such that $\E R'$ is a cyclic $R'$-module.
Then $\E R'=xR'$ for some $x\in\E^\reg$.
In particular, $R\subset y\E\subset R'$ for $y=x^{-1}\in Q_R^\reg$.
\end{lem}

\begin{proof}
By hypothesis, $R'$ is semilocal (see \cite[Exc.~9.3]{Mat89}), and $\E R'=zR'$ for some $z\in Q_R^\reg$.
Then $z^{-1}\E R'=R'$ implies the existence of a $w\in R'^*\cap z^{-1}\E$ (see the proof of \cite[Hilfssatz~2]{Jag77}), and $x:=zw$ satisfies the requirements.
\end{proof}

\subsection{Valuation rings}\label{80}

To deal with rings with zero-divisors, we need a general notion of valuation (ring), sometimes called a \emph{Manis} or \emph{pseudo-valuation (ring)} (see \cite{KV04,Mat73,CDK94}).
In the case of one-dimensional Cohen--Macaulay rings, only discrete valuation rings arise (see \S~\ref{77}).


\begin{dfn}\
\begin{enumerate}[label=(\alph*), ref=\alph*]

\item A ring $R$ is said to have a \emph{large Jacobson radical} if every prime ideal of $R$ containing the Jacobson radical of $R$ is a maximal ideal (see \cite[Ch.~I, (1.9) Prop.]{KV04}).

\item A ring $R$ is called \emph{Marot} if every regular ideal, or equivalently regular fractional ideal, $\E$ of $R$ is generated by $\E^\reg$.

\end{enumerate}
\end{dfn}


\begin{rmk}\label{19}\
\begin{enumerate}[label=(\alph*), ref=\alph*]

\item\label{19a} Any semilocal ring has a large Jacobson radical (see \cite[Ch.~I, (1.11) Rem.~(2)]{KV04}).

\item\label{19b} If $Q_R$ has a large Jacobson radical, then $R$ is a Marot ring (see \cite[Ch.~I, (1.12) Prop.]{KV04}).
In particular, this holds by \eqref{19a} if $R$ is reduced Noetherian.

\end{enumerate}
\end{rmk}


Let $Q$ be a ring with $Q^\reg=Q^*$ having a large Jacobson radical.


\begin{dfn}
A \emph{valuation ring} of $Q$ is a subring $V \subsetneq Q$ such that the set $Q \setminus V$ is multiplicatively closed.
For any ring $R\subset V$, we call $V$ a \emph{valuation ring over $R$}.
If $R\subset Q$ is a subring with $Q_R=Q$, we denote by $\V_R$ the set of all valuation rings of $Q$ over $R$.
\end{dfn}


\begin{rmk}\label{8}
Let $V$ be a valuation ring of $Q$.
\begin{enumerate}[label=(\alph*), ref=\alph*]
\item\label{8c} Then $V$ is integrally closed in $Q_V=Q$ (see \cite[Ch.~I, (2.1) Lem.]{KV04}).
\item\label{8a} There is a unique regular maximal ideal $\mm_V$ of $V$.
In particular, $V^\reg \setminus V^* \subset \mm_V$ (see \cite[Ch.~I, (2.2) Thm.]{KV04}).
\item\label{8b} Each $\E \in \RFI_V^*$ is principal (see \cite[Ch.~II, (2.2) Prop.~(2) and Ch.~I, (2.4) Prop.~(2)]{KV04}).
\end{enumerate}
\end{rmk}


Let $V$ be a valuation ring of $Q$.
Then the group $\RFI_V^\ast$ is totally ordered by reverse inclusion (see \cite[Ch.~I, (2.2) Thm.]{KV04}).
The \emph{infinite prime ideal} of $V$
\[
I_V := V : Q =\bigcap_{\E\in\RFI_V^\ast}\E \in \Spec(V)\cap\Spec(Q)
\]
is the intersection of all regular (principal) fractional ideals of $V$ (see \cite[Ch.~I, (2.4) Prop.~(3)(a)]{KV04}).
We include $\RFI_V^\ast$ into the totally ordered monoid 
\[
\RFI_{V,\infty}^\ast:=\RFI_V^\ast \cup \left\{ I_V \right\}.
\]
For $\E,\F\in\RFI_{V,\infty}^\ast$ we have $\E\F = I_V$ if $\left\{ \E, \F \right\} \not\subset \RFI_V^\ast$, and $\E < I_V$ for all $\E \in \RFI_V^\ast$.

For $x \in Q$, we denote by $\mu_V(x)$ the intersection of all regular $V$-submodules of $Q$ containing $x$.
If $x \in Q \setminus I_V$, then $\mu_V(x)\in \RFI_V^\ast$ (see \cite[Ch.~I, (2.4) Prop.~(3)(b)]{KV04}).
In particular, $\mu_V(x)=xV$ if $x \in Q^\reg$, and $\mu_V(x)=I_V$ if and only if $x \in I_V$.
This yields a map
\[
\mu_V \colon Q \to \RFI_{V,\infty}^\ast
\]
satisfying (see \cite[Ch.~I, (2.13) Prop.]{KV04})
\begin{equation}\label{9}
\mu_V ( x y ) = \mu_V ( x )\mu_V ( y ),\quad
\mu_V ( x + y ) \geq \min\{\mu_V(x),\mu_V(y)\}
\end{equation}
for any $x,y \in Q$, where equality holds if $\mu_V ( x ) \ne \mu_V ( y )$.
We can write
\begin{equation}\label{124a}
V = \left\{ x \in Q \mid \mu_V ( x ) \geq V \right\}
\end{equation}
with regular maximal ideal
\begin{equation}\label{124b}
\mm_V = \left\{ x \in Q \mid \mu_V ( x ) > V \right\}\supset I_V
\end{equation}
and units (see Remark~\ref{8}.\eqref{8a})
\begin{equation}\label{124c}
V^*=\left\{x\in Q^\reg\mid\mu_V(x)=V\right\}=(V\setminus\mm_V)^\reg.
\end{equation}


\begin{dfn}\label{10}
A valuation ring $V$ of $Q$ with regular maximal ideal $\mm_V$ is called a \emph{discrete valuation ring} if $\mm_V\in\RFI_V^\ast$ is the only regular prime ideal of $V$ (see \cite[Ch.~I, (2.16) Def.]{KV04}).
\end{dfn}
 

Let $V$ be a discrete valuation ring of $Q$. 
Then 
\begin{equation}\label{145}
\mm_V = \min\{\E\in\RFI_V^\ast\mid\E>V\}\in\RFI_V^*,
\end{equation}
$\bigcup_{k\in\ZZ}\mm_V^k=Q$, and $\bigcap_{k\in\ZZ}\mm_V^k=I_V$ (see \cite[Ch.~I, (2.15) Prop.]{KV04}). 
Therefore, there is a (unique) order preserving group isomorphism  
\begin{align}\label{11}
\phi_V\colon\RFI_V^* &\to\ZZ,\\
\nonumber\E & \mapsto\max\{k\in\ZZ\mid\mm_V^k\leq\E\},\\
\nonumber\mm_V^k & \mapsfrom k.
\end{align}
In fact, for $\E\in\RFI_V^\ast$ and $k\in\ZZ$ maximal with $\mm_V^k \leq \E$, we have $V=\mm_V^k:\mm_V^k\le\E : \mm_V^k<\mm_V$, and hence $\E=\mm_V^k$ by \eqref{145}.
Embedding $\ZZ$ into the totally ordered monoid
\[
\ZI:=\ZZ\cup\left\{\infty \right\}
\]
and extending $\phi_V$ by setting $\phi_V(I_V):=\infty$ yields a commutative diagram
\begin{equation}\label{12}
\xymatrix{
Q\ar@{>>}[d]_{\mu_V} \ar@{>>}[rd]^{\nu_V} \\
\RFI_{V,\infty}^\ast \ar[r]_{\cong}^{\phi_V} &\ZI.}
\end{equation}


\begin{dfn}\label{13}
A \emph{discrete valuation} of $Q$ is a map $\nu\colon Q\onto\ZI$ satisfying
\begin{equation*}
\nu (xy)=\nu(x)+\nu(y),\quad
\nu(x+y)\geq\min\{\nu(x),\nu(y)\},
\end{equation*}
for any $x,y \in Q$ (see \eqref{9}).
We refer to $\nu(x)\in\ZI$ as the \emph{value} of $x\in Q$ with respect to $\nu$.
The subring $V_\nu=\left\{x\in Q\mid\nu(x)\geq0\right\}$ of $Q$ is called the \emph{valuation ring} of $\nu$.
\end{dfn}


The valuation $\nu_V$ associated as above to a discrete valuation ring $V$ of $Q$ is discrete, and its valuation ring is $V_{\nu_V}=V$.

\section{Value semigroups}\label{77}

We specialize our setup to a one-dimensional semilocal Cohen--Macaulay ring $R$.
In this section we introduce value semigroups and value semigroup ideals, decompose them into local contributions, and show their invariance under completion.

\subsection{Admissible rings}\label{78}

One-dimensional local integrally closed Cohen-Macaulay rings are discrete valuation domains (see \cite[Ch.~II, (2.5) Prop.]{KV04}).
In general, the totality $\V_R$ of valuation rings of $Q_R$ over $R$ is described in the following theorem.
This provides the foundation for the definition and investigation of value semigroup ideals.


\begin{thm}\label{14}
Let $R$ be a one-dimensional semilocal Cohen--Macaulay ring.
\begin{enumerate}[label=(\alph*), ref=\alph*]
\item\label{14a} The set $\V_R$ is finite and non-empty, and it contains discrete valuation rings only.
\item\label{14b} We have $\Max(Q_R)=\{I_V\mid V\in\V_R\}$, and for any $I\in\Max(Q_R)$ there is a bijection
\begin{align*}
\{V\in\V_R\mid I_V=I\} &\to\V_{R/(I\cap R)},\\
V &\mapsto V/I,
\end{align*}
where $Q_{R/(I\cap R)}=Q_R/I$.
\item\label{14c} The integral closure of $R$ in $Q_R$ can be written as $\ol R =\bigcap\V_R$.
\item\label{14e} Any regular ideal of $\ol{R}$ is principal.
\item\label{14d} There is a bijection
\begin{align*}
\Max(\ol R) &\to \V_R, \\
\nn &\mapsto((\ol R\setminus\nn)^\reg)^{-1}\ol {R}, \\
\mm_V\cap\ol R &\mapsfrom V.
\end{align*}
In particular, $\ol R/(\mm_V\cap\ol R)=V/\mm_V$ and $\mm_V\cap R\in\Max(R)$.
\end{enumerate}
\end{thm}

\begin{proof}
See \cite[Ch.~II, (2.11) Thm.]{KV04} and use Lying Over for the particular claim of \eqref{14d}.
\end{proof}


By equation~\eqref{124c} and Theorem~\ref{14}.\eqref{14b} and \eqref{14c},
\begin{align}\label{123}
\ol R^*&=\{x\in Q_R\mid\nu(x) = 0\},\\
\nonumber R^*&=\ol{R}^*\cap R = \{x\in R\mid\nu(x)=0\}.
\end{align}

By Theorem~\ref{14}.\eqref{14e}, we have $\RFI_{\ol R} = \RFI_{\ol R}^*$, and there is a group isomorphism 
\begin{align}\label{131}
\psi=\psi_R\colon \RFI_{\ol{R}} &\to \prod_{V\in\V_R} \RFI_V^*, \\
\nonumber\E &\mapsto ( \E V)_{V\in \V_R},\\
\nonumber\bigcap_{V\in\V_R}\E_V &\mapsfrom(\E_V)_{V\in \V_R}.
\end{align}
In fact, writing $\E=t\ol R$ for some $t\in Q_R^\reg$,
\[
\bigcap_{V\in\V_R}\E V=\bigcap_{V\in\V_R}tV=t\bigcap_{V\in\V_R}V=t\ol R=\E
\]
by Theorem~\ref{14}.\eqref{14e}, and $\psi$ is injective.
Diagram~\eqref{12} taken component-wise with
\[
\phi=\phi_R:=\prod_{V\in\V_R}\phi_V
\]
gives rise to a commutative diagram
\begin{equation}\label{15}
\xymatrixcolsep{3pc}\xymatrix{
&Q_R^\reg \ar@{>>}[dl]\ar@{>>}[d]_-{\mu} \ar@{>>}[rd]^{\nu} \\
\RFI_{\ol{R}} \ar[r]_-{\cong}^-{\psi} &\displaystyle{\prod_{V\in \V_R} \RFI_V^*} \ar[r]_-{\cong}^-{\phi} &\ZZ^{\V_R}.
}
\end{equation}
Then surjectivity of $\nu$, and hence of $\psi$, follows from the approximation theorem for discrete valuations (see \cite[Ch.~I, (2.20) Thm.~(3)]{KV04}) which can be proved using Theorem~\ref{14}.\eqref{14d} and the Chinese remainder theorem.
The isomorphisms $\psi$ and $\phi$ preserve the partial orders on $\RFI_{\ol{R}}$ and $\prod_{V\in \V_R} \RFI_V^*$ by reverse inclusion and the natural partial order on $\ZZ^{\V_R}$.


\begin{dfn}\label{16}
Let $R$ be a one-dimensional semilocal Cohen--Macaulay ring, and let $\V_R$ be the set of (discrete) valuation rings of $Q_R$ over $R$ (see Theorem~\ref{14}.\eqref{14a}) with corresponding valuations
\[
\nu=\nu_R:=(\nu_V)_{V\in\V_R}\colon Q_R \to \ZZ_\infty^{\V_R}.
\]
\begin{enumerate}[label=(\alph*), ref=\alph*]

\item\label{16a} To each $\E\in\RFI_R$ we associate its \emph{value semigroup ideal}
\[
\Gamma_\E := \nu ( \E^\reg ) \subset \ZZ^{\V_R}.
\]
If $\E = R$, then the monoid $\Gamma_R$ is called the \emph{value semigroup} of $R$.

\item\label{16b} The value semigroup $\Gamma_R$ is called \emph{local} if $0$ is the only element of $\Gamma_R$ with a zero component in $\ZZ^{\V_R}$.

\item\label{16c}
We define a decreasing filtration $\Q^\bullet$ on $Q_R$ by 
\[
\Q^\alpha:=\{x\in Q_R\mid \nu(x)\ge\alpha\}
\]
for $\alpha\in\ZZ^{\V_R}$.
By $\E^\bullet:=\E\cap\Q^\bullet$ we denote the induced filtration on an $R$-submodule $\E$ of $Q_R$.

\end{enumerate}
\end{dfn}


\begin{lem}\label{127}
Let $R$ be a one-dimensional semilocal Cohen--Macaulay ring.
Then

\begin{enumerate}[label=(\alph*), ref=\alph*]

\item\label{127a} $\Q^\alpha=(\phi\circ\psi)^{-1}(\alpha)=\bigcap_{V\in\V_R}\mm_V^{\alpha_V}\in\RFI_{\ol R}$ for any $\alpha\in\ZZ^{\V_R}$,

\item\label{127b} with $\E$ also $\E^\alpha$ is a (regular) fractional ideal of $R$ for all $\alpha\in\ZZ^{\V_R}$,

\item\label{127c} $\Q^{\nu(x)}=x\ol R$ for any $x\in Q_R^\reg$ and, in particular, $\Q^0=\ol R$, and 

\item\label{127d} $\Gamma_{\Q^\alpha}=\alpha+\NN^{\V_R}$ for any $\alpha\in\ZZ^{\V_R}$ and, in particular, $\Gamma_{\ol R}=\NN^{\V_R}$.

\end{enumerate}
\end{lem}

\begin{proof}\
\begin{asparaenum}[(a)]

\item By definition of $\mu_V$, the first equality is due to isomorphism~\eqref{11} and diagram~\eqref{12}.
Isomorphisms~\eqref{11} and \eqref{131} yield the second equality.

\item Let $\E$ be a fractional ideal of $R$ and $\alpha\in\ZZ^{\V_R}$.
Then $\E^\alpha$ is an $R$-module by \eqref{127a}, and $r\E\subseteq R$ for some $r\in R^{\reg}$.
Thus, $r\E^\alpha\subseteq r\E\subseteq R$ and $\E^\alpha$ is a fractional ideal of $R$.
If $\E\in\RFI_R$, then there is an $x\in \E^{\reg}$.
By surjectivity of $\nu$ in diagram~\eqref{15} and equation~\eqref{124a}, there is a $y=z\frac yz\in(R^\beta)^\reg$ for arbitrarily large $\beta\in\ZZ^{\V_R}$.
Then $xy\in(\E^\alpha)^\reg$ for $\beta\ge\alpha-\nu(x)$, and hence $\E^\alpha\in\RFI_R$.

\item The particular claim is due to part~\eqref{127a} and Theorem~\ref{14}.\eqref{14c}.
The general claim then follows immediately by writing $y\in\Q^{\nu(x)}$ as $y=x\frac yx$ since $\nu(\frac yx)\ge0$.

\item The particular claim follows from surjectivity of $\nu$ in diagram~\eqref{15}, Theorem~\ref{14}.\eqref{14c}, and equation~\eqref{124a}.
Again by surjectivity of $\nu$, $\alpha=\nu(x)$ for some $x\in Q_R^\reg$.
Then the general claim follows using part~\eqref{127c}.
\qedhere

\end{asparaenum}
\end{proof}


The following result was stated without proof in \cite[(1.1.1)]{Del88} and \cite[\S2]{BDF00}.


\begin{prp}\label{100}
Let $R$ be a one-dimensional semilocal Cohen--Macaulay ring with value semigroup $\Gamma_R$.
Then $R$ is local if and only if $\Gamma_R$ is local.
\end{prp}

\begin{proof}
Suppose first that $R$ is local with maximal ideal $\mm$.
Then Theorem~\ref{14}.\eqref{14d} and equation~\eqref{124b} imply
\[
\mm\subset \bigcap_{V \in V_R} \mm_V=\bigcap_{V\in\V_R}\{x \in Q_R \mid \nu_V(x)>0\}.
\]
The statement follows with equation~\eqref{123}.

Suppose now that $\Gamma_R$ is local, and set $\mm:=\{x\in R\mid \nu(x)>0\}$.
By equation~\eqref{123}, any proper ideal of $R$ is contained in $\mm$.
Moreover, $\mm$ is obviously closed under multiplication by elements of $R$.
We show that $\nu(x)$ has no zero component for all $x\in\mm$. 
This implies that $\mm$ is also closed under addition, and hence an ideal.

For this, assume that there is an $x\in\mm$ such that $\nu_{V_1}(x) = 0$ for some $V_1 \in \V_R$.
Then $x \in R \setminus R^\reg\subset \bigcup_{V\in\V_R} I_V$ by hypothesis on $\Gamma_R$ and Theorem~\ref{14}.\eqref{14b}.
Thus, there is a $V_1 \ne V_2 \in \V_R$ such that $x \in I_{V_2}$.

By hypothesis on $R$, there is a $y \in R^\reg \setminus R^*$.
Then $\nu(y)\in\Gamma_R\setminus\{0\}$, and hence $\nu_V(y)>0$ for all $V\in\V_R$ by assumption on $\Gamma_R$.
After replacing $y$ by a suitable power, we may assume that $\nu_V(x) \ne \nu_V(y)$ for all $V\in \V_R$.
Then $\nu(x+y)=\min\{\nu(x),\nu(y)\}\in\ZZ^{\V_R}$.
Thus, $x + y \in R^\reg$ again since $R \setminus R^\reg \subset \bigcup_{V\in\V_R} I_V$, and hence $\nu (x+y)\in \Gamma_R$.

Therefore, by assumption on $\Gamma_R$, $\nu_{V_1} (x+y) = \nu_{V_1} (x) = 0$ yields $\nu (x+y) =0$, and thus $\nu_{V_2} (y) = \nu_{V_2}(x+y)=0$ contradicts the choice of $y$.
\end{proof}


In the following we show that, under suitable hypotheses, value semigroups $E=\Gamma_\E$ of fractional ideals $\E$ of $R$ satisfy certain axioms used to define the notion of good semigroup ideals in \S\ref{86}.


\begin{dfn}\label{88}
Let $R$ be a one-dimensional semilocal Cohen--Macaulay ring.
\begin{enumerate}[label=(\alph*), ref=\alph*]

\item\label{88d} We call $R$ \emph{analytically reduced} if $\wh R$ is reduced or, equivalently, $\wh{R_\mm}$ is reduced for all $\mm\in\Max(R)$ (see Lemma~\ref{5}.\eqref{5e}). 
 
\item\label{88a} The ring $R$ is called \emph{residually rational} if $\ol{R}/\nn=R/\nn\cap R$ for all $\nn\in\Max(\ol{R})$ or, equivalently, $V/\mm_V=R/\mm_V\cap R$ for all $V\in\V_R$ (see Theorem~\ref{14}.\eqref{14d}).

\item\label{88b} We say that $R$ has \emph{large residue fields} if $\lvert R/\mm \rvert \ge \lvert \V_{R_\mm} \rvert$ for all $\mm \in \Max(R)$.

\item\label{88c} We call $R$ \emph{admissible} if it is analytically reduced and residually rational with large residue fields.

\end{enumerate}
\end{dfn}


\begin{dfn}\label{17}
Let $\ol{S}$ be a partially ordered monoid, isomorphic to $\NN^I$ with its natural partial order, where $I$ is a finite set.
We consider the following properties of a subset $E$ of the \emph{group of differences} $D_{\ol S}\cong \ZZ^I$ of $\ol S$ (see \cite[\S1]{Del88} and \cite[\S2]{DAn97}).
\begin{enumerate}[label={(E\arabic*)}, ref=E\arabic*]
\setcounter{enumi}{-1}
\item\label{E0} There exists an $\alpha\in D_{\ol S}$ such that $\alpha+\ol S\subset E$.
\item\label{E1} If $\alpha,\beta \in E$, then $\min\{\alpha,\beta\} := (\min\{\alpha_i,\beta_i\})_{i\in I}\in E$.
\item\label{E2} For any $\alpha,\beta\in E$ and $j\in I$ with $\alpha_j=\beta_j$ there exists an $\varepsilon\in E$ such that $\varepsilon_j>\alpha_j=\beta_j$ and $\varepsilon_i\ge\min\{\alpha_i,\beta_i\}$ for all $i\in I\setminus\{j\}$ with equality if $\alpha_i\not=\beta_i$.
\end{enumerate}
We call $E$ \emph{good} if it satisfies \eqref{E0}, \eqref{E1}, and \eqref{E2}.
The \emph{difference} of $E,F\subset D_{\ol S}$ is
\[
E-F:=\{\alpha\in D_{\ol S}\mid\alpha+F\subset E\}.
\]
\end{dfn}


The following result shows that the isomorphism in Definition~\ref{17} is unique.


\begin{lem}\label{114}
Any group automorphism $\varphi$ of $\ZZ^s$ preserving the partial order is defined by a permutation of the standard basis.
\end{lem}

\begin{proof}
Let $\varphi$ be an automorphism of $\ZZ^s$ preserving the partial order.
Then $(\varphi(\ee_i))_{i\in \{1,\dots,s\}}$ is a basis of $\ZZ^s$, and hence $0<\ee_j=\sum_{i=1}^s\lambda_i\varphi(\ee_i)=\varphi(\sum_{i=1}^s\lambda_i\ee_i)$ for some $\lambda_i\in\ZZ$.
Since $\varphi$ is order preserving, this implies $\lambda_i\in\NN$ for all $i$.
For the $k$th component we have
\[
\sum_{i=1}^s \lambda_i \varphi (\ee_i)_k = (\ee_j)_k =
\begin{cases}
1 &\text{if } k = j, \\
0 &\text{otherwise.}
\end{cases}
\]
Since $\varphi$ is order preserving, we have $\varphi(\ee_i)>0$.
This yields $\ee_j = \varphi (\ee_i)$ for some $i$.
\end{proof}


\begin{lem}\label{110}
Let $R$ be a one-dimensional analytically reduced semilocal Cohen--Macaulay ring, and let $\E\in\RFI_R$.
Then $\ol R$ is a finite $R$-module, and hence $\RFI_{\ol R}\subset\RFI_R$.
In particular, $\C_\E\in\RFI_R\cap\RFI_{\ol R}$ and $\C_\E=x\ol R=\Q^{\nu(x)}$ for some $x\in Q^\reg$ with $\nu(x)+\NN^{\V_R}\subset\Gamma_\E$.
\end{lem}

\begin{proof}
If $R$ is analytically reduced, then $R$ is reduced (see Lemma~\ref{5}.\eqref{5a}). 
Hence, $Q_R$ localizes (see \cite[Cor.~2.1.13]{HS06}), and $\ol{R_\mm}$ is a finite $R_\mm$-module (see \cite[Ch.~II, (3.22) Thm.]{KV04}).
As integral closure localizes (see \cite[Prop.~2.1.6]{HS06}), we have $\ol{R_\mm}=\ol R_\mm$, and it follows that $\ol R$ is a finite $R$-module and hence $\RFI_{\ol R}\subset\RFI_R$.
The particular claim follows by Theorem~\ref{14}.\eqref{14e} and Lemma~\ref{127}.\eqref{127c} and \eqref{127d}.
\end{proof}


In the following, we collect results of D'Anna (see \cite{DAn97}) and provide a detailed proof.


\begin{prp}\label{18}
Let $R$ be a one-dimensional semilocal Cohen--Macaulay ring, and let $\E\in\RFI_R$.
\begin{enumerate}[label=(\alph*), ref=\alph*]
\item\label{18a} We have $\Gamma_\E+\Gamma_R\subset\Gamma_\E$.
\item\label{18b} If $R$ is analytically reduced, then $\Gamma_\E$ satisfies \eqref{E0} with $\ol S=\Gamma_{\ol R}$ and $I=\V_R$.
\item\label{18c} If $R$ is local and analytically reduced with large residue field, then $\Gamma_\E$ satisfies \eqref{E1}.
\item\label{18d} If $R$ is local and residually rational, then $\Gamma_\E$ satisfies \eqref{E2}.
\end{enumerate}
In particular, if $R$ is local admissible, then $\Gamma_\E$ is good (see Definition~\ref{17}).
\end{prp}

\begin{proof}\pushQED{\qed}\
\begin{asparaenum}[(a)]

\item Since $\E$ is an $R$-module and $Q_R^\reg=Q_R^*$ a group, $R^\reg\E^\reg\subset\E^\reg$.
Then the claim follows from $\nu$ in diagram~\eqref{15} being a group homomorphism.

\item Recall that $\ol S=\NN^I$ with $I=\V_R$ by Lemma~\ref{127}.\eqref{127d}, and $I$ is finite by Theorem~\ref{14}.\eqref{14a}.
By Lemma~\ref{110}, there is an $x \in Q^\reg$ such that
\[
\nu ( x ) + \ol S = \nu ( x \ol{R}^\reg ) = \nu (\C_\E^\reg) \subset \nu ( \E^\reg ) = \Gamma_\E.
\]

\item Let $x,y\in\E^\reg$ with $\nu(x)=\alpha$ and $\nu(y)=\beta$. 
By Lemma~\ref{110} and Theorem~\ref{14}.\eqref{14e} and \eqref{14d}, Lemma~\ref{7} applies to $R':=\ol R$.
We may thus assume that $\ideal{x,y}_{\ol R}=z\ol R$ for some $z\in\ideal{x,y}_R^\reg\subset\E^\reg$. 
Then $\nu(z)\ge\min\{\nu(x),\nu(y)\}\ge\nu(z)$ by Lemma~\ref{127}.\eqref{127d}, and hence $\min\{\nu(x),\nu(y)\}=\nu(z)\in\Gamma_\E$.

\item Denote by $\mm$ the maximal ideal of $R$.
Let $\alpha,\beta\in\Gamma_\E$ and $W\in\V_R$ such that $\alpha_W=\beta_W$.
Choose $x,y\in\E^\reg$ such that $\nu(x)=\alpha$ and $\nu(y)=\beta$. 
Then $\nu_W(x/y)=\alpha_W-\beta_W=0$, and hence $x/y\in W\setminus\mm_{W}$ by equations~\eqref{124a} and \eqref{124b}.
By hypothesis, $V/\mm_V=R/\mm$ for all $V\in\V_R$.
Thus, $\ol{x/y}=\ol{u}$ in $W/\mm_W=R/\mm$ for some $u\in R\setminus\mm$.
In particular, $\nu_W(u-x/y)>0$ and $\nu(u)=0$ by equations~\eqref{124a} and \eqref{124b}.
Then $uy-x\in\E$ with $\nu_W(uy-x)>\nu_W(y)=\beta_W$ and $\nu_V(uy-x)\ge\min\{\alpha_V,\beta_V\}$ for all $V \in\V_R\setminus\{W\}$ with equality if $\alpha_V\not=\beta_V$.
This remains true after replacing $u$ by any element $u'\in u+\mm$.
It is left to show that, for some $u'$, $\nu_V(u'-x/y)<\infty$ for all $V \in \V_R$ with $\alpha_V=\beta_V$.
Since $R$ is Cohen--Macaulay, there is a $z\in\mm^\reg\subset\mm_W^\reg$, and hence $(\infty,\dots,\infty)>\nu(z^k)\ge k\cdot(1,\dots,1)$.
Then $u'=u+z^k$ satisfies the requirement if $k>\nu_V(u-x/y)<\infty$ for all $V\in\V_R$ with $\alpha_V=\beta_V$.\qedhere
\end{asparaenum}
\end{proof}


While the value semigroup operation preserves inclusions, it is not compatible with the expected counterparts of multiplication and colon operation on the semigroup side.


\begin{rmk}\label{4}
Let $R$ be a one-dimensional semilocal Cohen--Macaulay ring, and let $\E,\F\in\RFI_R$.
\begin{enumerate}[label=(\alph*),ref=\alph*]

\item\label{4a} If $\E\subset\F$, then $\Gamma_{\E}\subset\Gamma_{\F}$.

\item\label{4b} The inclusion $\Gamma_{\E\F}\supset\Gamma_\E+\Gamma_\F$ is not an equality in general (see Example~\ref{118}).

\item \label{4c} The inclusion $\Gamma_{\E:\F}\subset\Gamma_\E-\Gamma_\F$  is not an equality in general (see \cite[Exa.~3.3]{BDF00}).

\end{enumerate}
\end{rmk}

\begin{figure}[ht]
\begin{tikzpicture}[inner sep=1.5,scale=0.5]
\draw[->] (0,0) -- (0,6);
\draw[->] (0,0) -- (10,0);
	
\foreach \i in {0,...,9} \foreach \j in {0,...,5} \draw (\i,\j) node[shape=circle,draw,fill=white] {};
\draw (0,0) node[shape=circle,draw,fill=black] {};
\foreach \i in {0,...,6} \foreach \j in {0,...,4} \draw (3+\i,1+\j) node[shape=circle,draw,fill=black] {};
	
\draw (5,6) node {$\Gamma_R$};
\end{tikzpicture}\quad
\begin{tikzpicture}[inner sep=1.5,scale=0.5]
\draw[->] (0,0) -- (0,6);
\draw[->] (0,0) -- (10,0);
	
\foreach \i in {0,...,9} \foreach \j in {0,...,5} \draw (\i,\j) node[shape=circle,draw,fill=white] {};
\foreach \i in {0,...,4} \draw (2,1+\i) node[shape=circle,draw,fill=black] {};
\draw (3,1) node[shape=circle,draw,fill=black] {};
\foreach \i in {0,...,4} \foreach \j in {0,...,3} \draw (5+\i,2+\j) node[shape=circle,draw,fill=black] {};
	
\draw (5,6) node {$\Gamma_\E$};
\end{tikzpicture}\medskip\\\noindent
\begin{tikzpicture}[inner sep=1.5,scale=0.5]
\draw[->] (0,0) -- (0,6);
\draw[->] (0,0) -- (10,0);
	
\foreach \i in {0,...,9} \foreach \j in {0,...,5} \draw (\i,\j) node[shape=circle,draw,fill=white] {};
\draw (3,1) node[shape=circle,draw,fill=black] {};
\foreach \i in {0,...,5} \foreach \j in {0,...,3} \draw (4+\i,2+\j) node[shape=circle,draw,fill=black] {};
	
\draw (5,6) node {$\Gamma_\F$};
\end{tikzpicture}\quad
\begin{tikzpicture}[inner sep=1.5,scale=0.5]
\draw[->] (0,0) -- (0,6);
\draw[->] (0,0) -- (10,0);
	
\foreach \i in {0,...,9} \foreach \j in {0,...,5} \draw (\i,\j) node[shape=circle,draw,fill=white] {};
\draw (5,2) node[shape=circle,draw,fill=black] {};
\draw (6,2) node[shape=circle,draw,fill=black] {};
\foreach \i in {0,...,4} \foreach \j in {0,1,2} \draw (5+\i,3+\j) node[shape=circle,draw,fill=black] {};
	
\draw (5,6) node {$\Gamma_\E+\Gamma_\F$};
\end{tikzpicture}
\caption{The value semigroup (ideals) in Example~\ref{118}.}\label{117}
\end{figure}
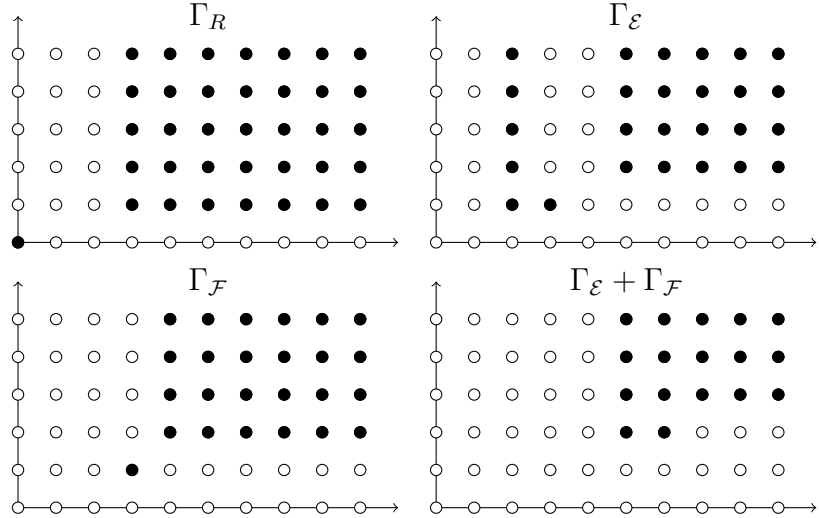

\begin{exa}\label{118}
Consider the admissible ring
\[
R := \CC [[ (-t_1^4,t_2),(-t_1^3,0),(0,t_2),(t_1^5,0) ]] \subset \CC [[t_1]] \times \CC [[t_2]] = \ol R
\] 
and the $R$-submodules of $Q_R$
\[
\E:= \langle (t_1^3,t_2),(t_1^2,0) \rangle_R,\quad\F := \langle (t_1^3,t_2), (t_1^4,0), (t_1^5,0) \rangle_R.
\]
Figure~\ref{117} shows that $R$ is local (see Proposition~\ref{100}), and that \eqref{E2} fails for $\Gamma_\E+\Gamma_\F$.
Thus, $\Gamma_{\E\F}\supsetneq\Gamma_\E+\Gamma_\F$ by Proposition~\ref{18}.
\end{exa}


\subsection{Compatibility with localization}\label{125}

Let $R$ be reduced.
Then $Q_R$, and hence $\ol R$, commutes with localization (see \cite[Cor.~2.1.13 and Prop.~2.1.6]{HS06}).


\begin{lem}\label{129}
Let $R$ be a reduced one-dimensional semilocal Cohen--Macaulay ring.
For any $\mm \in \Max(R)$ the localization map $\pi \colon Q_R \to (Q_R)_\mm = Q_{R_\mm}$ induces a bijection
\begin{align*}
\rho_\mm\colon\left\{ V \in \V_R \mid \mm_V \cap R = \mm \right\} &\to \V_{R_\mm}, \\
V &\mapsto V_\mm, \\
\pi^{-1} (W) &\mapsfrom W.
\end{align*}
In particular, $(\mm_V)_\mm=\mm_W$ if $V\mapsto W$.
\end{lem}

\begin{proof}
Let $\mm\in\Max(R)$ and $V \in \V_R$ with $\mm_V \cap R = \mm$, and hence $R \setminus \mm \subset V \setminus \mm_V$.
By exactness of localization, $(\mm_V)_\mm\subsetneq V_\mm$ contains a regular non-unit, and hence
\[
R_\mm \subset V_\mm \subsetneq (Q_R)_\mm=Q_{R_\mm}.
\]
An explicit calculation shows that $(Q_R)_\mm\setminus V_\mm$ is multiplicatively closed, and hence $V_\mm \in \V_{R_\mm}$.
Moreover, $V=\pi^{-1}(V_\mm)$ since $V\subsetneq Q_R$ is a maximal subring (see \cite[Ch.~I, (2.15) Prop.~(3)(d)]{KV04}).

Now, let $W\in\V_{R_\mm}$, and set $V:=\pi^{-1}(W)$.
Then $V_\mm=W\subsetneq Q_{R_\mm}$, and hence $R \subset V\subsetneq Q_R$.
With $Q_{R_\mm}\setminus W$ also $Q_R\setminus V$ is multiplicatively closed, and hence $V\in \V_R$.
Consider the commutative diagram of ring homomorphisms
\[
\xymatrix{
V\ar[r]^-\pi & W\\
R\ar[r]^-{\iota}\ar@{^(->}[u] & R_\mm.\ar@{^(->}[u]
}
\]
Using Theorem~\ref{14}.\eqref{14d}, $\pp:=\pi^{-1}(\mm_W)\in\Spec(V)$ satisfies 
\[
\pp\cap R=\iota^{-1}(\mm_W\cap R_\mm)=\iota^{-1}(\mm_{R_\mm})=\mm.
\]
In particular, with $\mm$ also $\pp$ is regular, and hence $\pp=\mm_V$ by Theorem~\ref{14}.\eqref{14a} (see Definition~\ref{10}).
\end{proof}

Let $R$ be a reduced one-dimensional semilocal Cohen--Macaulay ring. 
By Theorem~\ref{14}.\eqref{14d}, the sets $\left\{ V \in \V_R \mid \mm_V \cap R = \mm \right\}$, $\mm \in \Max(R)$, form a partition of $\V_R$.
By Lemma~\ref{129}, there is a bijection
\begin{align*}
\rho \colon \V_R &\to \bigsqcup_{\mm \in \Max(R)} \V_{R_\mm}, \\
V & \mapsto \rho_{\mm_V \cap R} (V) = V_{\mm_V \cap R},
\end{align*}
inducing an order preserving group homomorphism
\begin{align*}
\prod_{V \in \V_R} \RFI_V^* & \to\prod_{\mm \in \Max(R)} \prod_{W \in \V_{R_\mm}} \RFI_{W}^*, \\
\nonumber(\E_V)_{V \in \V_R} & \mapsto((\E_{\rho^{-1}(W)})_\mm)_{\mm\in\Max(R), W\in\V_{R_\mm}}.
\end{align*}
Since it maps $(\mm_V^{k_V})_{V \in \V_R}\mapsto(\mm_W^{k_{\rho^{-1}(W)}})_{\mm\in\Max(R), W\in\V_{R_\mm}}$, it is an isomorphism due to \eqref{11}.
Combined with diagram~\eqref{15} for $R$ and $R_\mm$ for $\mm\in\Max(R)$, it fits into a commutative diagram
\begin{equation*}
\xymatrix@C=3em{
Q_R^\reg \ar[ddd] \ar@{>>}[rd]\ar@/^1pc/@{>>}[rrrd]^-{\nu}\\
&\RFI_{\ol R} \ar[d]^-{\cong}_\xi \ar[r]_-{\cong}^-{\psi} &\displaystyle{\prod_{V \in \V_R} \RFI_V^*} \ar[d]^-{\cong} \ar[r]_-{\cong}^-{\phi} &\ZZ^{\V_R} \ar[d]^-{\cong} \\
&\displaystyle{\prod_{\mm \in \Max(R)} \RFI_{\ol{R_{\mm}}}} \ar[r]_-{\cong}^-{\prod_\mm\psi_{R_\mm}} &\displaystyle{\prod_{\mm \in \Max(R)} \prod_{W \in \V_{R_{\mm}}} \RFI_{W}^*} \ar[r]_-{\cong}^-{\prod_\mm\phi_{R_\mm}} &\displaystyle{\prod_{\mm \in \Max(R)}} \ZZ^{\V_{R_\mm}}\\
\displaystyle{\prod_{\mm \in \Max(R)} Q_{R_{\mm}}^\reg}\ar@{>>}[ru]\ar@/_1pc/@{>>}[rrru]_-{\prod_\mm\nu_{R_\mm}}}
\end{equation*}
where $\xi(\E):=(\E_\mm)_{\mm \in \Max(R)}$ with $\E_\mm\in\RFI_{\ol R_\mm}=\RFI_{\ol{R_\mm}}$ for any $\E\in\RFI_{\ol R}$ (see Lemma~\ref{2}.\eqref{2c}).
This implies
\begin{equation}\label{128}
\nu(x)=(\nu_{R_\mm}(x/1))_{\mm\in\Max(R)}
\end{equation}
for all $x \in Q_R^\reg$.
To ease notation, we identify the rightmost groups in the above diagram.


The first part of the following result was stated by Barucci, D'Anna, and Fr{\"o}berg (see \cite[\S~1.1]{BDF00}).


\begin{thm}\label{20}
Let $R$ be a one-dimensional reduced semilocal Cohen--Macaulay ring.
Then there is a decomposition into local value semigroups
\[
\Gamma_R = \prod_{\mm \in \Max(R)} \Gamma_{R_\mm},
\]
and any $\E\in\RFI_R$ decomposes as
\[
\Gamma_\E = \prod_{\mm \in \Max(R)} \Gamma_{\E_\mm}.
\]
\end{thm}

\begin{proof}
By Proposition~\ref{100}, $\Gamma_{R_\mm}$ is local for all $\mm \in \Max(R)$.
It remains to prove the second decomposition.

By equation~\eqref{128}, there is an inclusion $\Gamma_\E\subset\prod_{\mm \in \Max(R)} \Gamma_{\E_\mm}$.
Let now $\alpha=(\alpha_\mm)_{\mm \in \Max(R)} \in \prod_{\mm \in \Max(R)} \Gamma_{\E_\mm}$.
Then there are $x_\mm/y_\mm\in \E_\mm$, $\mm\in\Max(R)$, such that $\nu_{R_\mm} (x_\mm/y_\mm)=\alpha_\mm$ for every $\mm \in \Max(R)$.
By equation~\eqref{123}, we may clear denominators and assume that $y_\mm=1$ for every $\mm \in \Max(R)$.
By the Chinese remainder theorem, there is a $z_\mm\in(\bigcap \Max(R) \setminus\{\mm\})\setminus\mm$ for each $\mm \in \Max(R)$.
Then by equations~\eqref{124a} and \eqref{124b} and by Theorem~\ref{14}.\eqref{14d}, we have $\nu_{R_\mm}(z_\mm/1) = 0$ and $\nu_V(z_\mm/1)>0$ for all $V \in \V_{R_\nn}$ for every $\nn \in \Max(R)\setminus\{\mm\}$.
Let
\[
k_\mm > \max \{\nu_V (x_\nn/1) - \nu_V (x_\mm/1) \mid V \in \V_{R_\nn},\ \nn \in \Max(R)\setminus\{\mm\}\}.
\]
Then $z:=\sum_{\mm \in \Max(R)} x_\mm z_\mm^{k_\mm} \in \E$ with $\nu_{R_\mm} (z/1) = \alpha_\mm$ for any $\mm \in \Max(R)$.
Thus, $\nu(z)=\alpha$ by equation~\eqref{128}.
The claimed equality follows.
\end{proof}


\begin{cor}\label{87}
Let $R$ be a one-dimensional reduced semilocal Cohen--Macaulay ring with large residue fields, and let $\E \in \RFI_R$.
\begin{enumerate}[label=(\alph*), ref=\alph*]
\item\label{87a} If $R$ is analytically reduced, then $\Gamma_\E$ satisfies \eqref{E1}.
\item\label{87b} If $R$ is residually rational, then $\Gamma_\E$ satisfies \eqref{E2}.
\end{enumerate}
In particular, if $R$ is admissible, then $\Gamma_\E$ is good (see Definition~\ref{17}).
\end{cor}

\begin{proof}
Using Theorem~\ref{20}, this follows from Proposition~\ref{18}.
Note that to prove property \eqref{E2} for elements $\alpha,\beta\in\Gamma_\E$ which are different in all components in $\Gamma_{\E_\mm}$ for some $\mm \in \Max(R)$, we need to apply \eqref{E1} in $\Gamma_{\E_\mm}$.
\end{proof}

\subsection{Invariance under completion}

The invariance of value semigroup ideals with completion is due to D'Anna (see \cite[\S1]{DAn97}).
We give a proof including the semilocal case.


\begin{lem}\label{130}
With $R$ also $\wh R$ is a one-dimensional (semi)local Cohen--Macaulay ring.
\end{lem}

\begin{proof}
See Lemma~\ref{5}.\eqref{5e} and \cite[Cor.~2.1.8]{BH93}.
\end{proof}


\begin{thm}\label{21}
Let $R$ be a one-dimensional local Cohen--Macaulay ring.
Then there is a bijection
\begin{align*}
\sigma\colon\V_R &\to \V_{\wh R},\\
V &\mapsto V\wh R,\\
W\cap Q_R &\mapsfrom W.
\end{align*}
In particular, $\mm_V\wh R=\mm_W$ if $V\mapsto W$.
\end{thm}

\begin{proof}
See \cite[Ch.~II, (3.19) Thm.~(2)]{KV04} for the bijection.

With $\mm_V$ also $\mm_V\wh R$ is regular, and $R/\mm_R=\wh R/\mm_R\wh R$ implies $\wh R=R+\mm_R\wh R$ (see Lemma~\ref{5}.\eqref{5a} and \eqref{5b}).
Since $\mm_R\subset\mm_V$ and hence $V\mm_R\wh R\subset\mm_V\wh R$, it follows that $V\wh R/\mm_V\wh R=V/(\mm_V\wh R\cap V)=V/\mm_V$ (see Lemma~\ref{2}.\eqref{2d}).
The particular claim thus follows by uniqueness of $\mm_W$ (see Remark~\ref{8}.\eqref{8a}).
\end{proof}


\begin{cor}\label{22}
Let $R$ be a one-dimensional local Cohen--Macaulay ring. 
Then $\ol{\wh{R}}=\ol{R}\wh{R}$.
In particular, $\ol{\wh{R}}=\wh{\ol{R}}$ if $\ol R$ is finite over $R$.
\end{cor}

\begin{proof}
This follows from Lemma~\ref{130} and Theorems~\ref{14}.\eqref{14c} and \ref{21} using Lemmas~\ref{2}.\eqref{2d} and \ref{5}.\eqref{5a} (see \cite[Ch.~II, (3.19) Thm.~(3)]{KV04}).
Lemmas~\ref{2}.\eqref{2c} and \ref{5}.\eqref{5f} yield the particular claim.
\end{proof}


Let $R$ be a one-dimensional local Cohen--Macaulay ring. 
By Theorem~\ref{21}, there is an order preserving group homomorphism
\begin{align*}
\prod_{V\in \V_R} \RFI_V^* &\to \prod_{W\in \V_{\wh R}} \RFI_{W}^*\\
(\E_V)_{V\in\V_R} &\mapsto (\E_{\sigma^{-1}(W)}\wh R)_{W\in \V_{\wh R}}
\end{align*}
mapping $(\mm_V^{k_V})_{V \in \V_R}\mapsto(\mm_W^{k_{\sigma^{-1}(W)}})_{W\in\V_{\wh R}}$ which is an isomorphism due to \eqref{11}.
Combined with diagram~\eqref{15} for $R$ and $\wh R$ (see Lemma~\ref{130}), it fits into a commutative diagram
\[
\xymatrix@C=5em{
Q_R^\reg \ar[ddd]\ar@{>>}[rd]\ar@/^1pc/@{>>}[rrrd]^-{\nu}\\
 &\RFI_{\ol{R}} \ar[d]_-{\eta}^-{\cong} \ar[r]^-{\psi}_-{\cong} &\displaystyle{\prod_{V\in \V_R} \RFI_V^*} \ar[d]^-{\cong} \ar[r]^-{\phi}_-{\cong} &\ZZ^{\V_R} \ar[d]^-{\cong}\\
 &\RFI_{\ol{\wh{R}}} \ar[r]^-{\psi_{\wh R}}_-{\cong} &\displaystyle{\prod_{W\in \V_{\wh R}} \RFI_{W}^*} \ar[r]^-{\phi_{\wh R}}_-{\cong} &\ZZ^{\V_{\wh R}}\\
Q_{\wh{R}}^\reg\ar@{>>}[ru]\ar@/_1pc/@{>>}[rrru]_-{\nu_{\wh R}}
}
\]
where $\eta\colon\E\mapsto\E\wh R$ and $\eta^{-1}\colon\F\cap Q_R\mapsfrom\F$.
To ease notation, we identify the rightmost groups in the above diagram.


The following lemma relates value semigroup ideals to jumps in the filtration induced by $\Q^\bullet$ (see \cite[(4.3) Rem.]{CDK94}).


\begin{lem}\label{133}
Let $R$ be a one-dimensional analytically reduced local Cohen--Macaulay ring with large residue fields, and let $\E\in\RFI_R$.
\begin{enumerate}[label=(\alph*), ref=\alph*]

\item\label{133a} For any $\alpha\in\ZZ^{\V_R}$, $\alpha \in \Gamma_\E$ is equivalent to $\E^\alpha/\E^{\alpha+\ee_V} \ne 0$ for all $V \in \V_R$.

\item\label{133b} If $R$ is residually rational, then $\ell_R(\E^\alpha/\E^{\alpha+\ee_V})\le1$ for all $\alpha\in\ZZ^{\V_R}$ and $V\in \V_R$.

\end{enumerate}
\end{lem}

\begin{proof}\
\begin{asparaenum}[(a)]

\item By Remark~\ref{19}.\eqref{19b}, $R$ is a Marot ring, and by Lemma~\ref{127}.\eqref{127b}, $\E^\alpha\in\RFI_R$ is generated by regular elements.
Thus, $\E^\alpha/\E^{\alpha+\ee_V}\ne 0$ if and only if there is a $\beta\in\Gamma_\E$ with $\beta\ge\alpha$ and $\beta_V=\alpha_V$.
The claim follows since $\Gamma_\E$ satisfies property~\eqref{E1} by Proposition~\ref{18}.\eqref{18c}.

\item By diagram~\eqref{15}, $\alpha=\nu(x)$ for some $x\in Q_R^\reg$.
By Definition~\ref{16}.\eqref{16c}, there is an isomorphism
\[
\xymatrix{
\E^\alpha/\E^{\alpha+\ee_V}\subset \Q^\alpha/\Q^{\alpha+\ee_V} & \Q^0/\Q^{\ee_V}=\ol R/(\mm_{V}\cap\ol R)=V/\mm_{V}\ar[l]_-{\cdot x}^-{\cong}
}
\]
for every $V\in \V_R$, where the equalities are due to Lemma~\ref{127}.\eqref{127a} and \eqref{127c} and Theorem~\ref{14}.\eqref{14d}.
If $R$ is residually rational, then $V/\mm_{V}=R/\mm$, and the claim follows.
\qedhere

\end{asparaenum}
\end{proof}


\begin{thm}\label{24}
Let $R$ be a one-dimensional analytically reduced semilocal Cohen--Macaulay ring with large residue fields.
Then
\[
\Gamma_\E = \Gamma_{\wh{\E}}
\]
for any $\E \in \RFI_R$.
\end{thm}

\begin{proof}
By Lemma~\ref{5}.\eqref{5b} and \eqref{5e},
 \[
\wh{\E_\mm}
=\E\otimes_RR_\mm\otimes_{R_\mm}\wh{R_\mm}
=\E\otimes_R\wh{R_\mm}
=\E\otimes_R\wh R_{\wh\mm}
=\E\otimes_R\wh R\otimes_{\wh R}\wh R_{\wh\mm}
=\wh\E_{\wh\mm}
\]
for all $\mm\in\Max(R)$.
Therefore, Theorem~\ref{20} reduces the claim to the case where $R$ is local.

By Lemmas~\ref{2}.\eqref{2d}, \ref{5}.\eqref{5a}, and \ref{127}.\eqref{127a} and Theorem~\ref{21}, $\wh{\E^\alpha}=\wh\E^{\alpha}$ for all $\alpha\in\V_R$.
By Lemma~\ref{133}.\eqref{133a}, $\alpha \in \Gamma_\E$ is equivalent to $\E^\alpha/\E^{\alpha+\ee_V} \ne 0$ for all $V \in \V_R$.
This latter condition commutes with completion by Lemma~\ref{5}.\eqref{5a}.
\end{proof}

\section{Good semigroups}\label{86}

In this section, we consider semigroup ideals that satisfy the properties in Definition~\ref{17} which hold true in the case of value semigroup ideals (see Proposition~\ref{18} and Corollary~\ref{87}).
These semigroup ideals are called \emph{good} by Barucci, D'Anna, and Fr{\"o}berg~\cite{BDF00}.
As a combinatorial counterpart of the relative length of two fractional ideals, we describe the \emph{distance} of two good semigroup ideals.

\subsection{Axioms and properties}

Let $S$ be a cancellative commutative monoid. 
Then $S$ embeds into its (free abelian) group of differences $D_S$.
If $S$ is partially ordered, then $D_S$ carries a natural induced partial order.


\begin{dfn}\label{26}
Let $S$ be a partially ordered cancellative commutative monoid such that $\alpha \geq 0$ for all $\alpha \in S$.
Assume that $D_S$ is generated by a finite set $I$ such that there is an isomorphism $D_S\cong\ZZ^I$ which preserves the natural partial orders.
By Lemma~\ref{114}, $I$ is unique and serves to make sense of components of elements.

\begin{enumerate}[label=(\alph*), ref=\alph*]

\item\label{26a}
If $0$ is the only element of $S$ with a zero component, then 
we call $S$ \emph{local}.

\item\label{26b}
We call $S$ a \emph{good semigroup} if it satisfies properties~\eqref{E0}, \eqref{E1}, and \eqref{E2} with
\[
\ol S := \{ \alpha \in D_S \mid \alpha \geq 0 \} \cong \NN^I.
\]
\item\label{26c}
A \emph{semigroup ideal} of a good semigroup $S$ is a subset $\emptyset\ne E \subset D_S$ such that $E + S \subset E$ and $\alpha+E\subset S$ for some $\alpha\in S$.

\end{enumerate}

Let now $E$ and $F$ be semigroup ideals of a good semigroup $S$.

\begin{enumerate}[label=(\alph*), ref=\alph*]
\setcounter{enumi}{3}

\item\label{26d}
We call 
\[
C_E:=E-\ol S=\{ \alpha \in D_S \mid \alpha + \ol S \subset E \}
\]
the \emph{conductor (semigroup) ideal} of $E$ and set $C:=C_S$.

\item\label{26e}
If $E$ satisfies \eqref{E1}, then we denote by $\mu^E:=\min E$ its \emph{minimum} which exists due to Dickson's lemma~\cite{Dic13} and by $\gamma^E:=\mu^{C_E}$ its \emph{conductor}.
Note that
\[
C_E=\gamma^E+\ol S.
\]
We abbreviate $\tau^E:=\gamma^E-(1,\dots,1)$, $\gamma:=\gamma^S$ and $\tau:=\tau^S$.

\item\label{26f}
If $E$ satisfies \eqref{E1} and \eqref{E2}, then we call $E$ a \emph{good semigroup ideal} of $S$.
The set of good semigroup ideals of $S$ is denoted by $\GSI_S$.

\end{enumerate}
\end{dfn}


\begin{rmk}\label{103}\
\begin{enumerate}[label=(\alph*), ref=\alph*]

\item\label{103a} Any semigroup ideal $E$ of $S$ satisfies property~\eqref{E0} since $S$ does and $E+S\subset E$.

\item\label{103b} If $S\subset S'\subset\ol S$ are good semigroups, then $D_{S'}=D_S$ and hence $\ol{S'}=\ol S$.
It follows that $\GSI_{S'}\subset\GSI_{S}$ and, in particular, $S'\in\GSI_{S}$.

\item\label{103d}
For any semigroup ideal $E$ of $S$ satisfying \eqref{E1}, $\mu^E = 0$ is equivalent to $S \subset E \subset \ol S$.
In fact, if $\mu^E = 0$, then $S = 0 + S = \mu^E + S \subset E$, and $\alpha \geq \mu^E = 0$ for all $\alpha \in E$ implies $E \subset \ol S$.
Conversely, if $S \subset E \subset \ol S$, then $0 = \mu^S \geq \mu^E \geq \mu^{\ol S} = 0$.

\item\label{103c}
Let $R$ be an admissible (local) ring. 
Then $S=\Gamma_R$ satisfies property~\eqref{E0} with $\ol S=\Gamma_{\ol R}=\NN^{\V_R}$ by Proposition~\ref{18}.\eqref{18b} and Lemma~\ref{127}.\eqref{127d}.
It follows that $D_{\Gamma_R}=D_{\ol S}=\ZZ^{\V_R}$.
Then $S$ is a good (local) semigroup, and $\Gamma_\E\in \GSI_S$ for any $\E\in\RFI_R$ by Proposition~\ref{100} and Corollary~\ref{87}.

\end{enumerate}
\end{rmk}


We collect some trivial properties of the difference for future reference.


\begin{rmk}\label{102}
Let $S$ be a good semigroup, $\alpha\in D_S$, and $E,E',F,F'$ be semigroup ideals of $S$.
Then
\begin{enumerate}[label=(\alph*), ref=\alph*]
\item\label{102d} $E-S=E$,
\item\label{102a} $\GSI_S\to\GSI_S$, $E\mapsto\alpha+E$, is an inclusion preserving bijection,
\item\label{102b} $(\alpha+E)-F=\alpha+(E-F)=E-(-\alpha+F)$, and
\item\label{102c} $E-F' \subset E-F \subset E'-F$ if $E\subset E'$ and $F\subset F'$.
\end{enumerate}
\end{rmk}


Although $\GSI_S$ is neither a monoid nor closed under difference (see Remark~\ref{4}), there is at least the following positive result (see Lemma~\ref{110}).


\begin{lem}\label{29}
For any two semigroup ideals $E$ and $F$ of a good semigroup $S$, also $E-F$ is a semigroup ideal of $S$.
If $E$ satisfies \eqref{E1}, so does $E-F$, and $C_E \in \GSI_S \cap \GSI_{\ol S}$.
\end{lem}

\begin{proof}
Since $F$ is a semigroup ideal of $S$, we have $(E-F)+S+F=(E-F)+F\subset E$, and hence $(E-F)+S\subset E-F$.
Since $E$ is a semigroup ideal of $S$, there is an $\alpha \in D_S$ such that $\alpha+E\subset S$.
Then we have for any $\beta \in F$, $\alpha+\beta+(E-F)\subset \alpha+E\subset S$.
Thus, $E-F$ is a semigroup ideal of $S$.

Assume now that $E$ satisfies property~\eqref{E1}.
Then for any $\alpha,\beta\in E-F$ and $\delta \in F$, we have $\min \{\alpha,\beta\}+\delta = \min \{\alpha+\delta,\beta+\delta\} \in E$ since $\alpha+\delta,\beta+\delta\in E$.
Hence, $\min\{\alpha,\beta\}\in E-F$, and $E-F$ satisfies property~\eqref{E1}.

We have $C_E+\ol S+\ol S=(E-\ol S)+\ol S+\ol S=(E-\ol S)+\ol S\subset E$, and hence $C_E+\ol S \subset E-\ol S = C_E$.
Therefore, $C_E$ is a semigroup ideal of $\ol S$.
As just shown, it satisfies \eqref{E1}, and hence $\min\{\alpha,\beta\}+\ol S\subset C_E$ for any $\alpha,\beta\in C_E$.
It follows that $C_E$ satisfies \eqref{E2}.
\end{proof}


\begin{rmk}\label{141}
Let $M$ be a finite index set, and let $S_m$, $m\in M$, be good semigroups.
Then $S=\prod_{m\in M}S_m$ is a good semigroup, $D_S=\prod_{m \in M}D_{S_m}$, and $\ol S = \prod_{m \in M} \ol{S_m}$.
Let $E_m$ denote the image of $E\subset D_S$ under projection to $D_{S_m}$.

Let $E=\prod_{m \in M}E_m$ and $F=\prod_{m \in M}F_m$.
Then $E-F=\prod_{m \in M}(E_m-F_m)$.
If, for all $m\in M$, $E_m$ is a (good) semigroup ideal of $S_m$, then $E$ is a (good) semigroup ideal of $S$.
In particular, $\gamma^E = (\gamma^{E_m})_{m\in M}$ if the latter is defined.
\end{rmk}


The following result decomposes good semigroups and their good semigroup ideals into local components.


\begin{thm}\label{109}
Any good semigroup $S$ decomposes uniquely and compatibly with the partial orders as a finite direct product 
\[
S=\prod_{m\in M}S_m
\]
of good local semigroups $S_m$.
Any semigroup ideal $E$ of $S$ satisfying \eqref{E1} decomposes as (see Remark~\ref{141})
\[
E=\prod_{m\in M}E_m.
\]
In particular, if $E\in\GSI_S$, then $E_m\in\GSI_{S_m}$ for all $m\in M$.
\end{thm}

\begin{proof}
See \cite[Thm.~2.5, Rem.~2.6, Prop.~2.12]{BDF00}.
\end{proof}


\begin{rmk}\label{142}
As value semigroups and their ideals are special good semigroups and good semigroup ideals (see Corollary~\ref{87}), the decompositions in Theorem~\ref{20} are special cases of those in Theorem~\ref{109}.
\end{rmk}


The following objects were introduced by Delgado~\cite{Del87,Del88} for investigating the Gorenstein symmetry.
They detect equality in Lemma~\ref{133}.\eqref{133b} in the case where $E=\Gamma_\E$ (see \cite[(4.6) Rem.]{CDK94}).


\begin{dfn}\label{101}
Let $S$ be a good semigroup, $E$ a semigroup ideal of $S$, $\alpha \in D_S$, and $J \subset I$.
We set
\begin{enumerate}[label=(\alph*), ref=\alph*]
\item $\Delta_J (\alpha):=\{\beta\in D_S \mid\alpha_i=\beta_i \text{ for } i\in J \textup{ and } \alpha_j<\beta_j \text{ for } j\in I\setminus J\}$,
\item $\Delta_J^E (\alpha) := \Delta_J (\alpha) \cap E$,
\item $\Delta(\alpha):=\bigcup_{i\in I} \Delta_i (\alpha)$, where $\Delta_i(\alpha):=\Delta_{\{i\}}(\alpha)$, and
\item $\Delta^E (\alpha) := \Delta (\alpha) \cap E$.
\end{enumerate}
\end{dfn}


In the remainder of this subsection, we provide some technical tools for \S\ref{136}.
The following two lemmas were proved by Delgado in the case where $E=S$ (see \cite[Lem.~1.8 and Cor.~1.9]{Del88}).


\begin{lem}\label{31}
Let $S$ be a good semigroup and $E \in \GSI_S$.
Assume that there is a $\delta \in E$ and a $J\subset I$ such that $\delta_j \geq \gamma^E_j$ for all $j \in J$.
If $\alpha \in D_S$ with
\begin{align*}
 \alpha_j &\geq \gamma^E_j \text{ for all } j \in J, \\
 \alpha_i &= \delta_i \text{ for all } i \in I \setminus J,
\end{align*}
then $\alpha \in E$.
\end{lem}

\begin{proof}
Choose an $\varepsilon \in D_S$ such that
\begin{align*}
 \varepsilon_j &= \delta_j \text{ for all } j \in J, \\
 \varepsilon_i &> \max\{\gamma^E_i, \delta_i \} \text{ for all } i \in I \setminus J.
\end{align*}
In particular, $\varepsilon \geq \gamma^E$, and hence $\varepsilon \in E$.
By property~\eqref{E2} applied to $\delta$ and $\varepsilon$, we obtain for any $j \in J$ a $\delta' \in E$ with $\delta' \geq \delta + \ee_j$ and $\delta'_i=\delta_i$ for all $i\in I\setminus J$.
Therefore, we may assume that $\delta \geq \alpha$.

Choose an $\varepsilon \in D_S$ such that
\begin{align*}
 \varepsilon_j &= \alpha_j \text{ for all } j \in J, \\
 \varepsilon_i &> \max\{\gamma^E_i, \alpha_i\}\text{ for all } i \in I \setminus J.
\end{align*}
In particular, $\varepsilon \geq \gamma^E$, and hence $\varepsilon \in E$. 
Thus, $\alpha=\min\{\varepsilon,\delta\}\in E$ since $E$ satisfies \eqref{E1}.
\end{proof}


\begin{lem}\label{32}
Let $S$ be a good semigroup. 
Then $\Delta^E ( \tau^E ) = \emptyset$ for any $E \in \GSI_S$.
\end{lem}

\begin{proof}
Assume that $\Delta^E ( \tau^E ) \ne \emptyset$. 
Then there is an $i \in I$ with a $\delta \in \Delta_i^E ( \tau^E )$.
That is, $\delta_i = \gamma^E_i-1$ and $\delta_j \geq \gamma^E_j$ for all $j \in I \setminus \{i\}$.
Therefore, Lemma~\ref{31} implies $\gamma^E-\ee_i+\ol S \subset E$, contradicting the minimality of $\gamma^E$ in $C_E=E-\ol S$.
\end{proof}


\begin{lem}\label{33}
Let $E$ and $F$ be semigroup ideals of a good semigroup $S$ satisfying property~\eqref{E1}.
Then $\gamma^{E-F} = \gamma^E-\mu^F$.
\end{lem}

\begin{proof}
Note that $\gamma^{E-F}$ is defined since $E - F$ satisfies property~\eqref{E1} by Lemma~\ref{29}. 
Since $F - \mu^F \subset \ol S$ and $\gamma^E + \ol S \subset E$, we have 
\[
\gamma^E - \mu^F + \ol S+F\subset\gamma^E + \ol S \subset E, 
\]
and hence $\gamma^E - \mu^F \geq \gamma^{E-F}$.
Conversely, $\gamma^{E-F}+\mu^F\ge\gamma^E$ follows from
\[
\gamma^{E-F} + \mu^F+\ol S = \gamma^{E-F}+\mu^F-\mu^F+F+\ol S=\gamma^{E-F}+\ol S+F \subset E.\qedhere
\]
\end{proof}

\subsection{Distance and length}\label{84}

\begin{dfn}\label{113}
Let $S$ be a good semigroup, and let $E\subset D_S$ be a subset.

Two elements $\alpha,\beta\in E$ are called \emph{consecutive} in $E$ if $\alpha<\beta$ and $\alpha<\delta<\beta$ implies $\delta\not\in E$ for any $\delta\in D_S$.
A chain 
\begin{equation}\label{91}
\alpha=\alpha^{(0)}<\cdots<\alpha^{(n)}=\beta
\end{equation}
of elements $\alpha^{(i)}\in E$ is said to be \emph{saturated} of \emph{length} $n$ if $\alpha^{(i)}$ and $\alpha^{(i+1)}$ are consecutive in $E$ for all $i \in \{0,\dots, n-1\}$.
Let $E$ satisfy
\begin{enumerate}[label={(E\arabic*)}, ref=E\arabic*]
\setcounter{enumi}{3}
\item\label{E4} For any fixed $\alpha,\beta\in E$, any two saturated chains \eqref{91} in $E$ have the same length.
\end{enumerate}
Then the \emph{distance} of $\alpha$ and $\beta$ in $E$ with $\alpha\le\beta$ is the length
\[
d_E(\alpha,\beta):=n
\]
of a saturated chain~\eqref{91}.
The \emph{distance} between two semigroup ideals $E \subset F$ of $S$ satisfying properties~\eqref{E1} and \eqref{E4} is then
\[
d(F\backslash E):=d_F(\mu^F,\gamma^E)-d_E(\mu^E,\gamma^E).
\]
\end{dfn}


\begin{prp}\label{89}
Let $S$ be a good semigroup. 
Then any $E \in \GSI_S$ satisfies property~\eqref{E4}.
\end{prp}

\begin{proof}
See \cite[2.3~Prop.]{DAn97}.
\end{proof}


\begin{rmk}\label{93}
Let $S$ be a good semigroup, and let $E\subset F$ be semigroup ideals of $S$ satisfying properties~\eqref{E1} and \eqref{E4}.

\begin{enumerate}[label=(\alph*), ref=\alph*]

\item\label{93c} $d_E$ is additive with respect to composition of chains.

\item\label{93d} $d_E(\alpha,\beta)\le d_F(\alpha,\beta)$ for all $\alpha,\beta\in E$ with $\alpha\le\beta$.

\item\label{93e} $d(F\backslash E)=d(\alpha+F\backslash\alpha+E)$ for all $\alpha\in D_S$.

\item\label{93a} In the situation of Theorem~\ref{109} (see \cite[Prop.~2.12.(iii)]{BDF00}) with $E,F\in\GSI_S$, we have
\[
d(F\backslash E)=\sum_{m\in M} d(F_m\backslash E_m).
\]

\item\label{93b} If $\varepsilon\ge\gamma^E$, then
\begin{align*}
d(F\backslash E)
&=d_F(\mu^F,\gamma^E)-d_E(\mu^E,\gamma^E) \\
&=d_F(\mu^F,\gamma^E)+d_F(\gamma^E,\varepsilon)-d_E(\mu^E,\gamma^E)-d_E(\gamma^E,\varepsilon) \\
&=d_F(\mu^F,\varepsilon)-d_E(\mu^E,\varepsilon)
\end{align*}
by \eqref{93c} and since $d_F(\gamma^E,\varepsilon)=d_E(\gamma^E,\varepsilon)$.

\end{enumerate}
\end{rmk}


In the following, we collect the main properties of the distance function $d$.
We begin with additivity.


\begin{lem}\label{38}
Let $S$ be a good semigroup, and let $E \subset F$ be two semigroup ideals of $S$ satisfying properties~\eqref{E1} and \eqref{E4}.
Then
\[
d(G\backslash E)=d(G\backslash F)+d(F\backslash E).
\]
\end{lem}

\begin{proof}
This can be seen using Remark~\ref{93}.\eqref{93b} (see \cite[2.7~Prop.]{DAn97}).
\end{proof}


The distance function detects equality as formulated by D'Anna (see \cite[2.8~Prop.]{DAn97}).
The proof of this fact is an immediate consequence of the following lemma.


\begin{lem}\label{34}
Let $E\subset F$ be two semigroup ideals of a good semigroup $S$, where $E\in\GSI_S$ and $F$ satisfies property~\eqref{E1}.
Let $\alpha \in F\backslash E$ be minimal.
Then any $\beta\in E$ maximal with $\beta<\alpha$ and $\beta'\in E$ minimal with $\alpha<\beta'$ are consecutive in $E$.
\end{lem}

\begin{proof}
Suppose that $\beta<\varepsilon<\beta'$ for some $\varepsilon\in E$.
By choice of $\beta$ and $\beta'$, $\alpha\not\le\varepsilon\not\le\alpha$, and hence $\min\{\alpha,\varepsilon\}<\alpha$.
By property~\eqref{E1} of $F$, $\min\{\alpha, \varepsilon\}\in F$, and hence $\min\{\alpha, \varepsilon\}\in E$ by minimality of $\alpha\in F\setminus E$, and $\beta=\min\{\alpha,\varepsilon\}$ by maximality of $\beta$.
In particular, $\beta_j=\varepsilon_j<\alpha_j\le\beta_j'$ for some $j\in I$.
Applying property~\eqref{E2} to $\beta,\varepsilon \in E$ yields an $\varepsilon'\in E$ with $\beta<\varepsilon'$, where $\beta_j<\varepsilon'_j$.
After replacing $\varepsilon'$ by $\min\{\varepsilon',\beta'\}\in E$, using property~\eqref{E1} of $E$, $\beta<\varepsilon'<\beta'$, and hence $\beta=\min\{\alpha,\varepsilon'\}$.
However, this contradicts $\beta_j<\alpha_j,\varepsilon'_j$.
\end{proof}


\begin{prp}\label{37}
Let $S$ be a good semigroup, and let $E,F\in\GSI_S$ with $E\subset F$.
Then $E=F$ if and only if $d(F\backslash E)=0$.
\end{prp}


\begin{proof}
For the non-trivial implication, assume that $d(F\backslash E)=0$ but $E\subsetneq F$.
In particular, $\mu^F=\mu^E$ by Remark~\ref{93}.\eqref{93c}.
Choose a minimal $\alpha \in F\backslash E$.
Then $\mu^E<\alpha<\gamma^E$.
In fact, assume that $\alpha\not\le\gamma^E$.
Then applying property~\eqref{E1} of $F$ to $\alpha$ and $\gamma^E$ yields a $\delta\in F$ with $\delta<\alpha,\gamma^E$, and hence $\delta \in E$ by minimality of $\alpha$.
However, Lemma~\ref{31} then implies that $\alpha\in E$, contradicting the assumption on $\alpha$.
By Lemma~\ref{34}, we have
\[
\mu^F=\mu^E\le\beta<\alpha<\beta'\le\gamma^E
\]
for some consecutive $\beta,\beta'\in E$.
By Proposition~\ref{89} and Remarks~\ref{93}.\eqref{93c} and \eqref{93d},
\begin{align*}
d_F(\mu^F,\gamma^E)
&=d_F(\mu^F,\beta)+d_F(\beta,\beta')+d_F(\beta',\gamma^E)\\
&>d_E(\mu^E,\beta)+d_E(\beta,\beta')+d_E(\beta',\gamma^E)
=d_E(\mu^E,\gamma^E),
\end{align*}
contradicting the hypothesis.
\end{proof}


Finally, we show that the distance function coincides with the relative length of fractional ideals when evaluated on their value semigroup ideals.


\begin{prp}\label{39}
Let $R$ be an admissible ring.
If $\E, \F \in \RFI_R$ such that $\E\subset\F$, then
\[
\ell_R(\F/\E)=d(\Gamma_\F\backslash \Gamma_\E).
\]
\end{prp}

\begin{proof}
See \cite[2.2~Prop.]{DAn97} for part of the following proof in the local case.
By Corollary~\ref{87}, $E:=\Gamma_\E$ and $F:=\Gamma_\F$ are good semigroup ideals of $\Gamma_R$, and hence by Corollary \ref{89}, they satisfy property~\eqref{E4}.

Let $\rr$ be the Jacobson radical of $R$.
By Theorem~\ref{14}.\eqref{14d}, $\rr\subset\bigcap_{V\in\V_R}\mm_V$ and hence $\nu(x)\ge(1,\dots,1)$ for all $x\in\rr$ by equation~\eqref{124b}.
By Lemma~\ref{110}, $\C_\E=\Q^{\varepsilon}$ for some $\varepsilon\in\ZZ^{\V_R}$ with $\varepsilon\ge\gamma^E$.
It follows that, for suffciently large $k\in\NN$, 
\[
\rr^k\F\subset \Q^{\mu^F+k\cdot(1,\dots,1)}\subset\Q^\varepsilon=\C_\E\subset\E.
\]
This turns $\F/\E$ into a module over the product ring (see Lemma~\ref{5}.\eqref{5e})
\[
R/\rr^k=\prod_{\mm\in\Max(R)}R_\mm/\mm^k.
\]
It follows that $\F/\E=\prod_{\mm\in\Max(R)}(\F/\E)_\mm$, and hence 
\[
\ell_R(\F/\E)=\sum_{\mm\in\Max(R)}\ell_{R_\mm}(\F_\mm/\E_\mm).
\]
By Theorem~\ref{20} and Remark~\ref{93}.\eqref{93a}, we may therefore assume that $R$ is local.
By Lemma~\ref{133}, then $\ell_R(\E^\alpha/\E^{\alpha+\ee_V})\le1$ with equality for all $V\in \V_R$ if and only if $\alpha\in E$.

Let $\alpha,\beta\in E$ be consecutive in $E$.
Then $d_E(\alpha,\beta)=1$ by definition.
For any $\delta\in\ZZ^{\V_R}$ with $\alpha<\delta<\beta$, $\delta\not\in E$ and hence $\ell_R(\E^\delta/\E^{\delta+\ee_V})=0$ for some $V\in \V_R$.
If $\delta_W=\beta_W$ for some $W\in \V_R$, then $\E^\beta/\E^{\beta+\ee_W}\subset\E^\delta/\E^{\delta+\ee_W}$ and hence $\ell_R(\E^\delta/\E^{\delta+\ee_W})\ge\ell_R(\E^\beta/\E^{\beta+\ee_W})=1$ since $\beta\in E$.
Thus, $\delta_V<\beta_V$ and hence $\ell_R(\E^\alpha/\E^\beta)=1$ by additivity of length.

By additivity of length and distance, it follows that
\[
d_E(\mu^E,\varepsilon)=\ell_R(\E^{\mu^E}/\E^{\varepsilon})=\ell_R(\E/\E^{\varepsilon})
\]
for any $\varepsilon\ge\gamma^E$, and hence (see Remark~\ref{93}.\eqref{93b})
\begin{align*}
d(F\backslash E)&=d_F(\mu^F,\varepsilon)-d_E(\mu^E,\varepsilon) \\
&=\ell_R(\F/\F^{\varepsilon})-\ell_R(\E/\E^{\varepsilon})\\
&=\ell_R(\F/\E^{\varepsilon})-\ell_R(\E/\E^{\varepsilon})
=\ell_R(\F/\E).\qedhere
\end{align*}
\end{proof}


As a consequence, the value semigroup ideal detects equality of regular fractional ideals (see \cite[2.5~Cor.]{DAn97}).


\begin{cor}\label{30}
Let $R$ be an admissible ring, and let $\E,\F \in \RFI_R$ such that $\E\subset\F$. 
Then $\E = \F$ if and only if $\Gamma_\E = \Gamma_\F$.
\end{cor}

\begin{proof}
See Propositions~\ref{37} and \ref{39}.
\end{proof}

\section{Dualities}\label{81}

This section is devoted to duality and contains our main results.
After a review of canonical ideals, we develop a combinatoral duality on the good semigroup ideals of any good semigroup.
We show that it mirrors the duality by canonical ideals by taking values.

\subsection{Cohen--Macaulay duality}\label{79}

Let $R$ be a one-dimensional Cohen--Macaulay ring.
In the following we recall some basics of canonical ideals.
We begin with the definition (see \cite[Def.~2.4]{HK71}).


\begin{dfn}\label{40}
Let $R$ be a one-dimensional Cohen--Macaulay ring.
A regular fractional ideal $\K \in \RFI_R$ is said to be a \emph{canonical (fractional) ideal} of $R$ if, for all $\E \in \RFI_R$,
\[
\E=\K:(\K:\E)
\]
or, equivalently, $\E=\Hom_R(\Hom_R(\E,\K),\K)$ (see Lemma~\eqref{2}.\eqref{2a}).
In particular, $R=\K:\K$.
\end{dfn}


Dualizing with a canonical ideal preserves relative length of regular fractional ideals.


\begin{lem}\label{73}
Let $R$ be a one-dimensional Cohen--Macaulay ring, $\K$ a canonical ideal of $R$ and $\E,\F\in \RFI_R$ with $\E\subset\F$. 
Then
\[
\ell_R(\K:\E/\K:\F)=\ell_R(\F/\E).
\]
\end{lem}

\begin{proof}
See \cite[Rem.~2.5.(c)]{HK71}.
\end{proof}


Being a canonical ideal is a local property in the following sense.


\begin{lem}\label{72}
Let $R$ be a one-dimensional Cohen--Macaulay ring and $\K\in\RFI_R$.
Then $\K$ is a canonical ideal of $R$ if and only if $\K_\mm=\K R_\mm\in\RFI_{R_\mm}$ is a canonical ideal of $R_\mm$ for all $\mm\in\Max(R)$.
\end{lem}

\begin{proof}
This follows from Lemma~\ref{2} (see \cite[Lem.~2.6]{HK71}).
\end{proof}


\begin{rmk}\label{135}
Let $R$ be a one-dimensional Cohen--Macaulay ring, and let $\K$ be a canonical ideal of $R$.
For $\mm\in\Max(R)$, $\K_\mm$ is then of type $1$ by Lemma~\ref{73}.
In fact, if $R=(R,\mm)$ is local, then the type of $\K$ equals the length of
\[
\Ext_R^1(R/\mm,\K)\cong (\K:\mm)/(\K:R).
\]
Therefore, canonical ideals are canonical modules (see Remark~\ref{134} and \cite[Prop.~3.3.13 and Def.~3.3.16]{BH93}).
\end{rmk}


Canonical ideals are unique up to projective factors.


\begin{prp}\label{41}
Let $R$ be a one-dimensional Cohen--Macaulay ring with a canonical ideal $\K$.
Then $\K'$ is a canonical ideal of $R$ if and only if $\K'=\E\K$ for some invertible ideal $\E$ of $R$. 
In the case where $R$ is semilocal, the latter condition becomes $\K'=x\K$ for some $x\in Q_R^\reg$.
\end{prp}

\begin{proof}
See \cite[Satz 2.8]{HK71} and \S\ref{83}.
\end{proof}


If $R$ is local, then by Lemmas~\ref{2}.\eqref{2c}, \ref{5}.\eqref{5a}, and \ref{6}, $R$ has a canonical ideal $\K$ if and only if its completion $\wh{R}$ has a canonical ideal $\wh\K$ (see also \cite[Lem.~2.10]{HK71}).
This latter existence can be further characterized as follows.


\begin{thm}\label{90}
A one-dimensional local Cohen--Macaulay ring $R$ has a canonical ideal if and only if $\wh R$ is generically Gorenstein.
In particular, any one-dimensional analytically reduced local ring has a canonical ideal.
\end{thm}

\begin{proof}
See \cite[Kor.~2.12, Satz 6.21]{HK71}.
\end{proof}


\begin{cor}\label{97}
Any one-dimensional analytically reduced local Cohen--Macaulay ring $R$ with large residue field has a canonical ideal $\K$ such that $R\subset\K\subset\ol R$. 
It is unique up to multiplication by $\ol R^*$ with unique value semigroup ideal.
\end{cor}

\begin{proof}
By Theorem~\ref{90}, there is a canonical ideal $\E$ of $R$.
By Lemma~\ref{110} and Theorem~\ref{14}.\eqref{14e} and \eqref{14d}, Lemma~\ref{7} applies to $R'=\ol R$.
It yields a $y\in Q_R^\reg$ such that $\K:=y\E$ satisfies the inclusion requirements, and hence $\K\ol R=\ol R$.
By Proposition~\ref{41}, the canonical ideals of $R$ are of the form $\K'=x\K$ with $x\in Q_R^\reg$.
If also $\K'$ satisfies the inclusions, then $x\ol R=x\K'\ol R=\K\ol R=\ol R$, and hence $x\in\ol R^*$.
By \eqref{123}, $\nu(x)=0$ and thus $\Gamma_{\K'}=\Gamma_\K$.
\end{proof}


Finally, canonical ideals propagate along finite ring extensions (see \cite[Thm.~3.3.7.(b)]{BH93}).


\begin{lem}\label{42}
Let $\varphi\colon R\to R'$ be a local homomorphism of one-dimensional local Cohen--Macaulay rings such that $R'$ is a finite $R$-module and $Q_R=Q_{R'}$.
If $\K_R$ is a canonical ideal of $R$, then $\K_R:R'$ is a canonical ideal of $R'$.
\end{lem}

\begin{proof}
This follows from Remark~\ref{108}.\eqref{108d} and Definition~\ref{40}.
\end{proof}

\subsection{Duality on good semigroups}\label{136}

Let $S$ be a good semigroup.
Motivated by a result of J\"ager in the irreducible case (see \cite[Hilfssatz~5]{Jag77}), D'Anna introduced the following object (see \cite[\S3]{DAn97}) to characterize canonical ideals in terms of their value semigroup ideal (see Theorem~\ref{43}).


\begin{dfn}\label{46}
For any good semigroup $S$, we call (see Definitions~\ref{26}.\eqref{26e} and \ref{101})
\[
K^0_S:= \left\{ \alpha\in D_S \mid \Delta^S ( \tau-\alpha )=\emptyset \right\}.
\]
\emph{the (normalized) canonical (semigroup) ideal} of $S$.
\end{dfn}


\begin{lem}\label{47}
Let $S$ be a good semigroup.
Then the set $K^0_S$ is a semigroup ideal of $S$ satisfying property~\eqref{E1} with minimum $\mu^{K^0_S} = 0$ and conductor $\gamma^{K^0_S} = \gamma$.
\end{lem}

\begin{proof}
See \cite[3.2~Prop.]{DAn97} and Lemma~\ref{32}.
\end{proof}


Our definition of a canonical semigroup ideal below relies on the inclusion relations of good semigroup ideals and avoids a fixed conductor.


\begin{dfn}\label{44}
Let $S$ be a good semigroup (see Definition~\ref{26}).
We call $K\in\GSI_S$ a \emph{canonical (semigroup) ideal} of $S$ if $K \subset E$ implies $K=E$ for all $E \in \GSI_S$ with $\gamma^K=\gamma^E$.
\end{dfn}


\begin{rmk}\label{105}
By Remark~\ref{102}.\eqref{102a}, with $K$ also $\alpha + K$ is a canonical ideal of $S$ for any $\alpha \in D_S$.
\end{rmk}


The following result was stated by Barucci, D'Anna, and Fr{\"o}berg in the case where $K=K_S^0$ (see \cite[Prop.~2.15]{BDF00}).


\begin{prp}\label{112}
Let $S=\prod_{m\in M}S_m$ be the decomposition of a good semigroup $S$ into good local semigroups $S_m$ (see Theorem~\ref{109}).
A good semigroup ideal $K \in \GSI_S$ is a canonical ideal of $S$ if and only if $K_m$ is a canonical ideal of $S_m$ for every $m \in M$.
\end{prp}

\begin{proof}
First note that $K_m \in \GSI_{S_m}$ for any $m \in M$ by Theorem~\ref{109}.
Suppose that $K$ is a canonical ideal of $S$.
Let $m \in M$, and assume that $K_m$ is not a canonical ideal of $S_m$.
Then there is an $E_m \in \GSI_{S_m}$ with $\gamma^{E_m} = \gamma^{K_m}$ and $K_m \subsetneq E_m$.
By Remark~\ref{141}, $E:=E_m\times\prod_{n \in M \setminus \{m\}} K_n \in \GSI_S$ with $\gamma^E = \gamma^K$ and $K \subsetneq E$, contradicting $K$ being a canonical ideal.

Suppose now that $K_m$ is a canonical ideal of $S_m$ for all $m \in M$. 
Let $E \in \GSI_S$ with $\gamma^E = \gamma^K$ and $E \subset K$.
By Theorem~\ref{109} and Remark~\ref{141}, $E_m \in \GSI_{S_m}$ with $\gamma^{E_m} = \gamma^{K_m}$ and $K_m \subset E_m$ for all $m \in M$.
Since $K_m$ is a canonical ideal, this implies that $K_m = E_m$ for every $m \in M$, and hence $E = K$. Thus, $K$ is a canonical ideal.
\end{proof}


Our aim in this subsection is to establish the following result on canonical semigroup ideals in analogy with the ring case.


\begin{thm}\label{45}
Any good semigroup $S$ has a canonical ideal.
Moreover, for any $K \in \GSI_S$ the following statements are equivalent:
\begin{enumerate}[label=(\roman*),ref=\roman*]
\item\label{45a} $K$ is a canonical ideal of $S$.
\item\label{45b} There is an $\alpha \in D_S$ such that $\alpha+K = K^0_S$.
\item\label{45c} For all $E \in \GSI_S$ we have $K-(K-E)=E$.
\end{enumerate}
If $K$ is a canonical ideal of $S$, then the following hold:
\begin{enumerate}[label=(\alph*), ref=\alph*]
\item\label{45h} If $S\subset K\subset\ol S$, then $K=K^0_S$.
\item\label{45f} If $E \in \GSI_S$, then $K - E \in \GSI_S$.
\item\label{45d} $K - K = S$.
\item\label{45e} If $S' \subset \ol S$ is a good semigroup with $S \subset S'$, then $K' = K - S'$ is a canonical ideal of $S'$.
\end{enumerate}
\end{thm}

\begin{proof}\pushQED{\qed}
By Proposition~\ref{50}, $K^0_S\in\GSI_S$, and hence \eqref{45b}$\implies$\eqref{45a} yields existence.
\begin{asparaenum}[(i)]
\item[\eqref{45a} $\Rightarrow$ \eqref{45b}] See Proposition~\ref{56}.
\item[\eqref{45b} $\Rightarrow$ \eqref{45c}] See Remark~\ref{102}.\eqref{102b} and Proposition~\ref{62}.
\item[\eqref{45c} $\Rightarrow$ \eqref{45a}] See Proposition~\ref{61}.
\end{asparaenum}
\begin{asparaenum}[(a)]
\item[\eqref{45h}] See \eqref{45b}, Remark~\ref{103}.\eqref{103d} and Lemma~\ref{47}.
\item[\eqref{45f}] See \eqref{45b}, Remark~\ref{102}.\eqref{102a} and \eqref{102b} and Proposition~\ref{50}.
\item[\eqref{45d}] Set $E:=S$ in \eqref{45c}.
\item[\eqref{45e}] See Corollary~\ref{58}.\qedhere
\end{asparaenum}
\end{proof}

\begin{figure}[ht]
\begin{tikzpicture}[inner sep=1.5,scale=0.5]
\draw[->] (0,0) -- (0,9);
\draw[->] (0,0) -- (11,0);

\foreach \i in {0,...,10} \foreach \j in {0,...,8} \draw (\i,\j) node[shape=circle,draw,fill=white] {};
\foreach \i in {0,1,2} \draw (3*\i,2*\i) node[shape=circle,draw,fill=black] {};
\draw (5,4) node[shape=circle,draw,fill=black] {};
\foreach \i in {0,1,2} \draw (5,6+\i) node[shape=circle,draw,fill=black] {};
\foreach \i in {0,1,2} \foreach \j in {0,1,2} \draw (8+\i,6+\j) node[shape=circle,draw,fill=black] {};

\draw (5.5,9) node {$S$};
\end{tikzpicture}\quad
\begin{tikzpicture}[inner sep=1.5,scale=0.5]
\draw[->] (0,0) -- (0,9);
\draw[->] (0,0) -- (11,0);

\foreach \i in {0,...,10} \foreach \j in {0,...,8} \draw (\i,\j) node[shape=circle,draw,fill=white] {};
\foreach \i in {0,...,8} \draw (0,0+\i) node[shape=circle,draw,fill=black] {};
\foreach \i in {1,3,4,...,10} \draw (0+\i,0) node[shape=circle,draw,fill=black] {};
\draw (1,1) node[shape=circle,draw,fill=black] {};
\foreach \i in {0,...,6} \draw (3,2+\i) node[shape=circle,draw,fill=black] {};
\foreach \i in {1,...,7} \draw (3+\i,2) node[shape=circle,draw,fill=black] {};
\draw (4,3) node[shape=circle,draw,fill=black] {};
\foreach \i in {0,...,4} \draw (5,4+\i) node[shape=circle,draw,fill=black] {};
\foreach \i in {1,...,5} \draw (5+\i,4) node[shape=circle,draw,fill=black] {};
\foreach \i in {0,...,3} \draw (6,5+\i) node[shape=circle,draw,fill=black] {};
\draw (7,5) node[shape=circle,draw,fill=black] {};
\foreach \i in {0,...,2} \foreach \j in {0,...,2} \draw (8+\i,6+\j) node[shape=circle,draw,fill=black] {};

\draw (5.5,9) node {$K_S^0$};
\end{tikzpicture} \medskip\\
\begin{tikzpicture}[inner sep=1.5,scale=0.5]
\draw[->] (0,0) -- (0,9);
\draw[->] (0,0) -- (11,0);

\foreach \i in {0,...,10} \foreach \j in {0,...,8} \draw (\i,\j) node[shape=circle,draw,fill=white] {};
\foreach \i in {0,1,2} \draw (1+\i,2) node[shape=circle,draw,fill=black] {};
\foreach \i in {0,1,2} \draw (4+\i,4) node[shape=circle,draw,fill=black] {};
\foreach \i in {0,...,4} \foreach \j in {0,...,3} \draw (6+\i,5+\j) node[shape=circle,draw,fill=black] {};

\draw (5.5,9) node {$E$};
\end{tikzpicture}\quad
\begin{tikzpicture}[inner sep=1.5,scale=0.5]
\draw[->] (0,0) -- (0,9);
\draw[->] (0,0) -- (11,0);

\foreach \i in {0,...,10} \foreach \j in {0,...,8} \draw (\i,\j) node[shape=circle,draw,fill=white] {};
\foreach \i in {0,...,6} \draw (4+\i,2) node[shape=circle,draw,fill=black] {};
\draw (4,3) node[shape=circle,draw,fill=black] {};
\foreach \i in {0,...,3} \foreach \j in {0,...,4} \draw (7+\i,4+\j) node[shape=circle,draw,fill=black] {};

\draw (5.5,9) node {$K_S^0-E$};
\end{tikzpicture} \medskip\\
\begin{tikzpicture}[inner sep=1.5,scale=0.5]
\draw[->] (0,0) -- (0,9);
\draw[->] (0,0) -- (11,0);

\foreach \i in {0,...,10} \foreach \j in {0,...,8} \draw (\i,\j) node[shape=circle,draw,fill=white] {};
\foreach \i in {0,1,2} \draw (1+\i,2) node[shape=circle,draw,fill=black] {};
\foreach \i in {0,...,6} \foreach \j in {0,...,4} \draw (4+\i,4+\j) node[shape=circle,draw,fill=black] {};

\draw (5.5,9) node {$K_S^0 - (K_S^0-E)$};
\end{tikzpicture}
\caption{A semigroup ideal $E$ satisfying property~\eqref{E1} but not \eqref{E2}, where $K_S^0-E \not\in \GSI_S$ and $E\subsetneq K_S^0 - (K_S^0-E)$.}\label{119}
\end{figure}
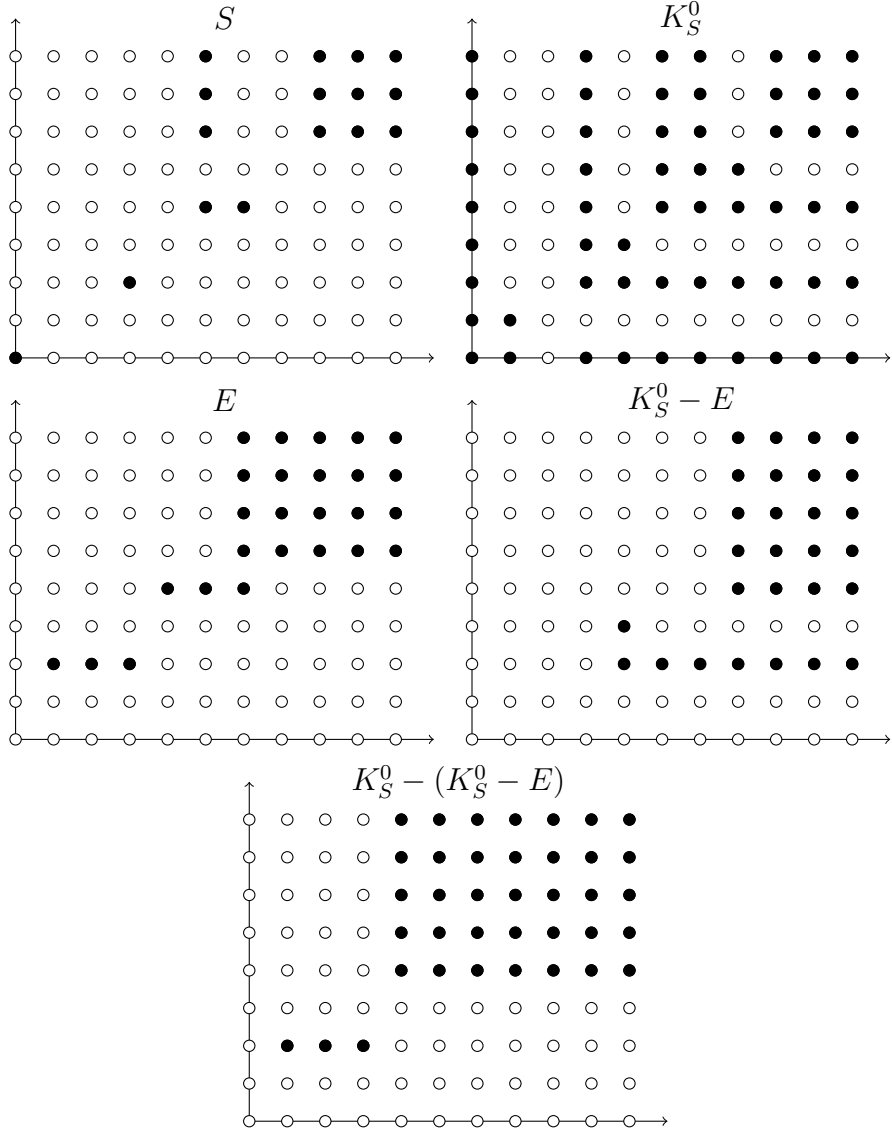

\begin{rmk}\label{120}
The assumption $E \in \GSI_S$ in Theorem~\ref{45}.\eqref{45c} and \eqref{45f} is necessary (see the example given by Figure~\ref{119}). 
\end{rmk}


We first approach Part~\eqref{45f} of Theorem~\ref{45} in the case where $K=K^0_S$.
To this end, we collect some properties of $K^0_S$.


\begin{lem}\label{48}\pushQED{\qed}
Let $S$ be a good semigroup.
Then the semigroup ideal $K^0_S$ of $S$ has the following properties:
\begin{enumerate}[label=(\alph*), ref=\alph*]
\item\label{48b} $\Delta^{K^0_S} (\tau) = \emptyset$.
\item\label{48a} If $E$ is a semigroup ideal of $S$, then 
\[K^0_S-E = \{\alpha \in D_S \mid \Delta^E(\tau-\alpha)=\emptyset\}.\]
\end{enumerate}
\end{lem}

\begin{proof}
This follows by calculation from Definitions~\ref{101} and \ref{46} (see \cite[3.3~Comp., 3.4~Lem.]{DAn97}).
\end{proof}


The proof of Theorem~\ref{45}.\eqref{45f} in the case where $K=K^0_S$ is achieved by the following proposition.
It shows, in particular, that $K^0_S$ is good.
D'Anna established a weaker statement, where \eqref{E2} is replaced by a certain property (E3) (see \cite[3.6~Thm.]{DAn97}).


\begin{prp}\label{50}
Let $S$ be a good semigroup.
Then $K^0_S-E\in\GSI_S$ for any $E \in \GSI_S$ and, in particular, $K^0_S\in \GSI_S$.
\end{prp}

\begin{figure}[ht]
\begin{tikzpicture}[tdplot_main_coords]

\draw[->] (0,0,0) -- (0,8,0); 
\draw (0,8,0) node[right] {$k_r$};
\draw[->] (0,0,0) -- (0,0,8); 
\draw (0,0,8) node[above] {$l_{r-1}$};
\draw[->] (0,0,0) -- (8,0,0); 

\draw (0,0,0) node[shape=circle,draw,fill=black,scale=0.5] {};
\draw (0,0,1) node[shape=circle,draw,fill=black,scale=0.5] {};

\draw[fill=black,opacity=0.3,line width=0] (1,1,0) -- (1,8,0) -- (8,8,0) -- (8,1,0) -- (1,1,0);
\foreach \i in {1,...,8} \foreach \j in {1,...,8} \draw (\j,\i,0) node[shape=circle,draw,fill=white,scale=0.5] {};
\draw (8,1,0) node[left, rotate=30] {$\Delta_{l_{r-1}}(\tau-\zeta^{(r)})$};

\draw[dashed] (4,0,1) -- (4,4,1);
\draw[dashed] (4,0,1) -- (6,0,1);
\draw[dashed] (4,0,1) -- (4,0,4);

\draw[fill=black,opacity=0.3,line width=0] (1,0,2) -- (8,0,2) -- (8,0,8) -- (1,0,8) -- (1,0,2);
\foreach \i in {2,...,8} \foreach \j in {1,...,8} \draw (\j,0,\i) node[shape=circle,draw,fill=white,scale=0.5] {};
\draw (8,0,5) node[left, rotate=30] {$\Delta_{k_r}(\tau-\zeta^{(r-1)})$};
\draw[opacity=0.3,line width=5pt] (1,0,1) -- (8,0,1);
\foreach \i in {1,...,8} \draw (\i,0,1) node[shape=circle,draw,fill=white,scale=0.5] {};
\draw (8,0,1) node[left, rotate=30] {$\Delta_{\{k_r,l_{r-1}\}}(\tau-\zeta^{(r-1)})$};

\draw[fill=black,opacity=0.3,line width=0] (0,1,2) -- (0,1,8) -- (0,8,8) -- (0,8,2) -- (0,1,2);
\foreach \i in {2,...,8} \foreach \j in {1,...,8} \draw (0,\j,\i) node[shape=circle,draw,fill=white,scale=0.5] {};
\draw[opacity=0.3,line width=5pt] (0,1,1) -- (0,8,1);
\foreach \i in {1,...,8} \draw (0,\i,1) node[shape=circle,draw,fill=white,scale=0.5] {};

\draw[fill=black,opacity=0.3,line width=0] (1,1,1) -- (1,8,1) -- (8,8,1) -- (8,1,1) -- (1,1,1);
\foreach \i in {1,...,8} \foreach \j in {1,...,8} \draw (\j,\i,1) node[shape=circle,draw,fill=white,scale=0.5] {};
\draw (8,1,1) node[left, rotate=30] {$\Delta_{l_{r-1}}(\tau-\zeta^{(r-1)})$};

\draw (4,4,1) node[shape=circle,draw,fill=black,scale=0.5] {};
\draw (4,4,1) node[above right] {$\delta^{(r-1)}$};

\draw (6,0,1) node[shape=circle,draw,fill=black,scale=0.5] {};
\draw (6,0,1) node[above] {$\delta^{(r)}$};

\draw (4,0,4) node[shape=circle,draw,fill=black,scale=0.5] {};
\draw (4,0,4) node[right] {$\varepsilon$};

\draw (0,0,0) node[above right] {$\tau-\zeta^{(r)}$};
\draw (0,0,1) node[above right] {$\tau-\zeta^{(r-1)}$};
\end{tikzpicture}
\caption{Induction step in the proof of Proposition~\ref{50} in the case where $I\setminus J=\{l_{r-1}\}$.}\label{116}
\end{figure}
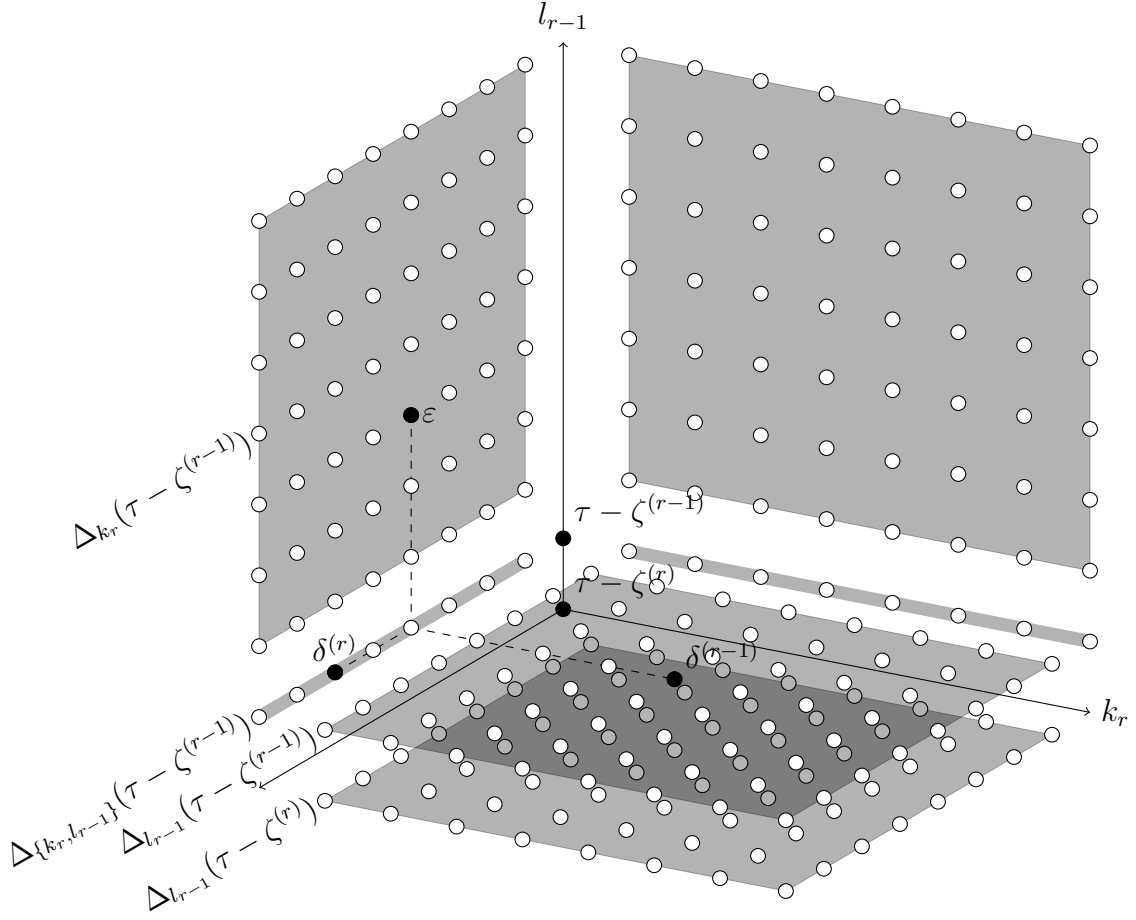

\begin{proof}
The idea of the following proof is illustrated in Figure~\ref{116}.

Suppose that $K^0_S - E \not\in \GSI_S$.
Since $K^0_S - E$ is a semigroup ideal of $S$ satisfying property~\eqref{E1} by Lemmas~\ref{29} and \ref{47}, it violates property~\eqref{E2}.
That is, there are $\alpha, \beta \in K^0_S - E$ with $\emptyset \ne J:=\{j\in I\mid \alpha_j\ne\beta_j\}\subset I$, $\zeta^{(0)}:=\min\{\alpha,\beta\}\in K^0_S-E$, and $l_0\in I \setminus J$ such that $\zeta\not\in K^0_S - E$ whenever $\zeta_{l_0} > \zeta_{l_0}^{(0)}$, $\zeta_i \ge \zeta_i^{(0)}$ for all $i \in I$, and $\zeta_j = \zeta_j^{(0)}$ for all $j \in J$.
In particular, any choice of a sequence $l_1,l_2,l_3,\ldots\in I\setminus J$ yields 
\begin{equation}\label{126}
\zeta^{(0)}\in K^0_S-E,\quad\zeta^{(r)}:=\zeta^{(r-1)}+\ee_{l_{r-1}}\not\in K^0_S-E.
\end{equation}
By Lemma~\ref{48}.\eqref{48a}, this means that $\Delta^E (\tau-\zeta^{(0)}) = \emptyset$, and, for all $r\ge1$, $\Delta_i^E(\tau-\zeta^{(r)})\ne\emptyset$ for some $i\in I$.
In order to construct a sequence of indices in $I\setminus J$ as above, we proceed by induction on $r$.
In each step we show additionally that $\Delta_j^E ( \tau-\zeta^{(r)} ) = \emptyset$ for all $j \in J$, and we choose an $l_r \in I \setminus J$ and a
\begin{equation}\label{122}
\delta^{(r)}\in\Delta_{l_r}^E(\tau-\zeta^{(r)}).
\end{equation}

Assume that this was done for $r-1$, and suppose that there is a $j \in J$ such that $\Delta_j^E(\tau-\zeta^{(r)})\ne\emptyset$.
Then $j\ne l_{r-1}$, and there is a
\begin{equation}\label{121}
\delta\in\Delta_j^E(\tau-\zeta^{(r)})
=\Delta_j^E(\tau-\zeta^{(r-1)})\sqcup\Delta_{\{j,l_{r-1}\}}^E(\tau-\zeta^{(r-1)})
=\Delta_{\{j,l_{r-1}\}}^E(\tau-\zeta^{(r-1)}),
\end{equation}
where the first equality holds by \eqref{126} and the second by the induction hypothesis.
We deduce a contradiction with different arguments for $r=1$ and $r\ge2$, respectively.

First consider the case $r=1$. 
Since $\beta \in K^0_S-E$ and $\delta\in E$, we get $\delta + \beta \in K^0_S$. 
Since $j\in J$, we may assume that $\beta_j>\zeta_j^{(0)}$, and we have $\beta_{l_0}=\zeta^{(0)}_{l_0}$ by choice of $l_0$.
By \eqref{121}, $\delta+\zeta^{(0)}\in\Delta_{\{j,l_0\}}(\tau)$ which implies that $\delta+\beta\in\Delta_{l_0}^{K^0_S}(\tau)$, contradicting Lemma~\ref{48}.\eqref{48b}.

Assume now that $r\ge2$.
By \eqref{122} and \eqref{121} and since $j\ne l_{r-1}$, 
\[
\delta^{(r-1)}_{l_{r-1}}=\tau_{l_{r-1}}-\zeta^{(r-1)}_{l_{r-1}}=\delta_{l_{r-1}},\quad
\delta^{(r-1)}_j>\tau_j-\zeta^{(r-1)}_j=\delta_j.
\] 
Then property~\eqref{E2} applied to $\delta^{(r-1)},\delta\in E$ yields an $\varepsilon \in E$ with $\varepsilon\geq\min\{\delta^{(r-1)},\delta\}\ge\tau-\zeta^{(r-1)}$, $\varepsilon_{l_{r-1}}>\delta_{l_{r-1}}$, and $\varepsilon_j=\delta_j$.
It follows that $\varepsilon \in \Delta^E_j(\tau- \zeta^{(r-1)})$, contradicting the induction hypothesis. 
 
Finally, choose an $r > \sum_{i \in I \setminus J} \lvert \tau_i -\zeta^{(1)}_i - \mu^E_i \rvert$.
Then $\delta^{(r)}_{l_r} = \tau_{l_r} - \zeta_{l_r}^{(r)} <\mu^E_{l_r}$ by \eqref{122}, contradicting the minimality of $\mu^E$.
It follows that $K^0_S-E\in\GSI_S$ as claimed.
With $E=S$, the particular claim follows by Remark~\ref{102}.\eqref{102d} and Lemma~\ref{47}.
\end{proof}


\tdplotsetmaincoords{70}{120}


We can now relate our canonical ideals (see Definition~\ref{44}) to D'Anna's normalized one (see Definition~\ref{46}).


\begin{prp}\label{56}
Let $S$ be a good semigroup, and let $K \in \GSI_S$. 
Then $K$ is a canonical ideal of $S$ if and only if $K = \alpha+K^0_S$ for some $\alpha \in D_S$.
\end{prp}

\begin{proof}
Using Remark~\ref{105}, it suffices to show that $K^0_S$ is the unique canonical ideal of $S$ with conductor $\gamma^{K^0_S}=\gamma$ (see Lemma~\ref{47}). 

To this end, we show that $E\subset K^0_S$ for any $E \in\GSI_S$ with $\gamma^E=\gamma$.
Since $K^0_S\in \GSI_S$ by Proposition~\ref{50}, this shows that $K^0_S$ is a canonical ideal.
Applied to any other canonical ideal $K$ of $S$ with conductor $\gamma^K=\gamma$, it gives $K = K^0_S$.

So let $E \in\GSI_S$ with $\gamma^E=\gamma$, and hence $\tau^E=\tau$.
Assume that there is a $\beta\in E\setminus K^0_S$.
Then there is a $\delta \in \Delta^S ( \tau - \beta )$ (see Definition~\ref{46}), and hence $\beta+\delta\in \Delta^E(\tau^E)$.
This contradicts Lemma~\ref{32}, and therefore $E\subset K^0_S$ as claimed.
\end{proof}


As a consequence we deduce the counterpart of Lemma~\ref{42} on the semigroup side (see Theorem~\ref{45}.\eqref{45e}). 


\begin{cor}\label{58}
Let $S \subset S'\subset\ol S$ be good semigroups.
If $K$ is a canonical ideal of $S$, then $K' = K - S'$ is a canonical ideal of $S'$.
\end{cor}

\begin{proof}
By Remark~\ref{103}.\eqref{103b}, $S' \in \GSI_S$, and, by Proposition~\ref{56}, $K = \alpha + K^0_S$ for some $\alpha \in D_S$.
Then by Lemma~\ref{48}.\eqref{48a},
\begin{align*}
K' &= ( \alpha + K^0_S ) - S' \\
   &= \alpha + ( K^0_S - S' ) \\
   &= \alpha + \{ \beta \in D_S \mid \Delta^{S'} ( \tau - \beta ) = \emptyset \} \\
   &= \alpha + \tau - \tau^{S'} + \{ \delta \in D_S \mid \Delta^{S'} ( \tau^{S'} - \delta ) = \emptyset \}.
\end{align*}
Thus, $K'$ is a canonical ideal of $S'$ by Proposition~\ref{56}.
\end{proof}


By the following two propositions, we establish an equivalent definition of canonical semigroup ideals (see Theorem~\ref{45}.\eqref{45c}) analogous to that of canonical ideals (see Definition~\ref{40}).


\begin{lem}\label{60}
Let $E$ and $F$ be semigroup ideals of a good semigroup $S$.
\begin{enumerate}[label=(\alph*), ref=\alph*]
\item\label{60a} There is an inclusion $E \subset F-(F-E)$.
\item\label{60b} If $E$ and $F$ satisfy property~\eqref{E1}, $F \subsetneq E$, and $\gamma^E = \gamma^F$, then $E \subsetneq F - ( F - E )$.
\end{enumerate}
\end{lem}

\begin{proof}\pushQED{\qedhere}\
\begin{asparaenum}[(a)]
\item This follows trivially from Definition~\ref{26}.\eqref{26d}.
\item Note that $F \subsetneq E$ forces $\mu^{F-E} > 0$.
Using Lemmas~\ref{33} and \ref{29}, we obtain
\[
\gamma^{F-( F - E )} = \gamma^F - \mu^{F-E} < \gamma^F = \gamma^E.
\]
Then the claim follows from \eqref{60a}. \qed
\end{asparaenum}
\end{proof}


\begin{prp}\label{61}
Let $S$ be a good semigroup, and let $K \in \GSI_S$ such that $K-(K-E)=E$ for all $E\in\GSI_S$. 
Then $K$ is a canonical ideal of $S$.
\end{prp}

\begin{proof}
Assume that $K$ is not a canonical ideal of $S$.
Then there is an $E\in\GSI_S$ with $\gamma^E = \gamma^K$ and $K\subsetneq E$ (see Definition~\ref{44}).
By Lemma~\ref{60}.\eqref{60b} and the hypothesis, this leads to the contradiction $E \subsetneq K - ( K - E )=E$.
\end{proof}


\begin{lem}\label{65}
Let $E$ be a semigroup ideal of a good semigroup $S$, and let $\alpha\in K^0_S-(K^0_S-E)$. 
If $\zeta \in D_S$ satisfies $\Delta^E(\tau-\zeta) = \emptyset$, then $\Delta^S(\tau-\zeta-\alpha)=\emptyset$.
Equivalently, if $\beta \in D_S$ satisfies $\Delta^S(\tau-\beta)\ne\emptyset$, then $\Delta^E(\tau-\beta+\alpha)\ne\emptyset$.
\end{lem}

\begin{proof}
Using Lemma~\ref{48}.\eqref{48a}, we compute
\begin{align*}
K^0_S - ( K^0_S - E ) &= \{ \alpha \in D_S \mid \alpha + ( K^0_S - E ) \subset K^0_S \} \\
&= \{ \alpha \in D_S \mid \alpha + \{ \zeta \in D_S \mid \Delta^E ( \tau - \zeta ) = \emptyset \} \subset K^0_S \} \\
&= \{ \alpha \in D_S \mid \forall\zeta\in D_S\colon\Delta^E(\tau-\zeta)=\emptyset\implies\Delta^S(\tau-\zeta-\alpha)=\emptyset\}.
\end{align*}
The equivalent formulation is obtained by setting $\zeta=\beta-\alpha\in D_S$.
\end{proof}


\begin{prp}\label{62}
Let $S$ be a good semigroup. 
Then $K^0_S - ( K^0_S - E ) = E$ for any $E \in \GSI_S$ and, in particular, $K^0_S - K^0_S = S$.
\end{prp}

\begin{proof}
By Lemma~\ref{60}.\eqref{60a}, there is a trivial inclusion
\[
E\subset K^0_S-(K^0_S-E)=:E'.
\]
By Lemmas~\ref{47} and \ref{29}, $E'$ is a semigroup ideal of $S$ satisfying condition~\eqref{E1}.
So in the case where $E\subsetneq E'$, there is a minimal $\alpha\in E'\setminus E$.
By property~\eqref{E1} of $E$, there is a $k \in I$ such that no $\varepsilon \in E$ satisfies $\varepsilon_k=\alpha_k$ and $\varepsilon_i\geq\alpha_i$ for all $i\in I\setminus\{k\}$.

We set $\beta:=\gamma-\ee_k\in D_S$, that is,
\begin{align*}
\beta_k&=\tau_k, \\
\beta_i&=\gamma_i \text{ for all } i \in I \setminus \{k\}.
\end{align*}
Then $0\in\Delta_k^S(\tau-\beta)$, and Lemma~\ref{65} yields a $\zeta\in\Delta_j^E( \tau - \beta + \alpha )$ for some $j\in I$.
That is, $\zeta \in E$ with
\begin{align*}
\zeta_j &= \tau_j  - \beta_j + \alpha_j, \\
\zeta_i &> \tau_i  - \beta_i + \alpha_i \text{ for all } i \in I \setminus \{j\}.
\end{align*}
We must have $j\ne k$ as otherwise $\varepsilon=\zeta$ would contradict the choice of $k$.
Thus,
\begin{align*}
\zeta_j &= \alpha_j - 1, \\
\zeta_k &> \alpha_k, \\
\zeta_i &\geq \alpha_i \text{ for all } i \in I \setminus \left\{ j,k \right\}.
\end{align*}
Since $\zeta\in E\subset E'$, by property~\eqref{E1} of $E'$ applied to $\zeta$ and $\alpha$, we find
\[
\alpha>\alpha-\ee_j =\min \{ \alpha, \zeta \}=:\alpha'\in E'.
\]
Property~\eqref{E2} of $E$ applied to $\alpha',\zeta\in E$ would yield an $\varepsilon\in E$ contradicting the choice of $k$.
Thus, $\alpha>\alpha' \in E'\setminus E$ contradicts the minimality of $\alpha$.
We conclude that $E=E'$.
With $E:=S$ the particular claim follows by Remark~\ref{102}.\eqref{102d} and Lemma~\ref{47}.
\end{proof}

\subsection{Relation of dualities}\label{137}

In this subsection, we put the Cohen--Macaulay duality in \S\ref{79} and the duality of good semigroup ideals in \S\ref{136} in relation.


We begin by extending the following result of D'Anna to semilocal rings.

\begin{thm}\label{43}
Let $R$ be an admissible local ring. 
Then a fractional ideal $\K$ of $R$ with $R\subset\K\subset\ol R$ is canonical if and only if $\Gamma_\K=K^0_{\Gamma_R}$ (see Definition~\ref{46}).
\end{thm}

\begin{proof}
See \cite[4.1~Thm.]{DAn97}.
\end{proof}


\begin{thm}\label{111}
Let $R$ be an admissible ring.
Then $\K\in\RFI_R$ is a canonical ideal of $R$ if and only if $\Gamma_\K$ is a canonical ideal of $\Gamma_R$.
\end{thm}

\begin{proof}
First assume that $R$ is local.
By Proposition~\ref{41} and Corollary~\ref{97}, $\K$ is a canonical ideal of $R$ if and only if there is an $x \in Q_R^\reg$ such that $x \K$ is a canonical ideal of $R$ with $R \subset x \K \subset \ol R$. 
By Theorem~\ref{43}, this is equivalent to $K_{\Gamma_R}^0 = \Gamma_{x \K} = \nu (x) + \Gamma_\K$. 
By Theorem~\ref{45}.\eqref{45a} $\Leftrightarrow$ \eqref{45b}, this is the case if and only if $\Gamma_\K$ is a canonical ideal of $\Gamma_R$.

Let now $R$ be semilocal.
By Lemma~\ref{72}, $\K$ is a canonical ideal of $R$ if and only if $\K_\mm$ is a canonical ideal of $R_\mm$ for every $\mm \in \Max(R)$. 
By the local case, this is equivalent to $(\Gamma_\K)_\mm = \Gamma_{\K_\mm}$ being a canonical ideal of $(\Gamma_R)_\mm = \Gamma_{R_\mm}$ for every $\mm \in\Max(R)$ (see Theorem~\ref{109} and Remark~\ref{142}). 
By Proposition~\ref{112} and Remark~\ref{142}, this is the case if and only if $\Gamma_\K$ is a canonical ideal of $\Gamma_R$.
\end{proof}


Next we show that taking values is compatible with the dualities of \S\ref{79} and \S\ref{136}.
We use the following result stated by Waldi in the case where $\E=R$ and $\F=\ol R$ (see \cite[Bem.~1.2.21]{Wal72}).


\begin{lem}\label{98}
Let $R$ be an admissible ring, $\E\in\RFI_R$ and $\F\in\RFI_{\ol R}$.
Set $E:=\Gamma_\E$ and $F:=\Gamma_\F$. 
Then $\E:\F=\Q^{\gamma^E-\mu^F}$ and hence $\Gamma_{\E:\F}=E-F$.
In particular, $\C_\E=\Q^{\gamma^E}$ and hence $\Gamma_{\C_\E}=C_E$.
\end{lem}

\begin{proof}
By Remark~\ref{108}.\eqref{108b}, Theorem~\ref{14}.\eqref{14e}, and Lemma~\ref{127}.\eqref{127c} and \eqref{127d}, it suffices to prove the particular claim.
By Lemma~\ref{110}, $\C_\E\subset \Q^{\gamma^E}$ (see Definition~\ref{26}.\eqref{26e}).
By Lemma~\ref{127}.\eqref{127d}, $\Gamma_{\Q^{\gamma^E}}\subset E$, and hence $\Gamma_{\E^{\gamma^E}}=\Gamma_{\Q^{\gamma^E}}=C_E$.
With Lemma~\ref{110} and Corollary~\ref{30}, it follows that $\Q^{\gamma^E}=\E^{\gamma^E}\subset\E$, and hence $\Q^{\gamma^E}\subset\C_\E$ since $\Q^{\gamma^E}$ is an $\ol R$-module.
\end{proof}


\begin{thm}\label{74}
Let $R$ be an admissible ring with canonical ideal $\K$.
Then
\begin{enumerate}[label=(\alph*), ref=\alph*]
\item\label{74a} $\Gamma_{\K:\F}=\Gamma_\K-\Gamma_\F$ for any $\F\in\RFI_R$ and
\item\label{74b} $d(\Gamma_\K-\Gamma_\E\backslash\Gamma_\K-\Gamma_\F)=d(\Gamma_\F\backslash \Gamma_\E)$ for any $\E,\F\in \RFI_R$ with $\E\subset\F$.
\end{enumerate}
\end{thm}

\begin{proof}\pushQED{\qed}
Set $S:=\Gamma_R$ and $K:=\Gamma_\K$.

By Theorem~\ref{20}, Lemmas~\ref{2}.\eqref{2b} and \ref{72}, and Remark~\ref{93}.\eqref{93a}, we may assume that $R$ is local.
By Remarks~\ref{108}.\eqref{108b}, \ref{102}.\eqref{102b} and \ref{93}.\eqref{93e}, Proposition~\ref{41}, Corollary~\ref{97}, and Theorem~\ref{43}, we may further assume that $K=K^0_S$.

We now prove both \eqref{74a} and \eqref{74b} simultaneously, setting $\E:=\F$ in the first case.
By Proposition~\ref{39} and Lemma~\ref{73},
\[
d(\Gamma_{\K:\E}\backslash \Gamma_{\K:\F})=\ell_R((\K:\E)/(\K:\F))=\ell_R(\F/\E)=d(\Gamma_\F\backslash \Gamma_\E)=:n.
\]
In particular, since $\C_\E\in\RFI_R$ by Lemma~\ref{110},
\begin{equation}\label{140}
d(\Gamma_{\K:\C_\E}\backslash\Gamma_{\K:\F})=\ell_R(\F/\C_\E)=:m+n.
\end{equation}
Choose a composition series in $\RFI_R$ (see \cite[Ch.~6]{AM69})
\[
\C_\E=\E_0\subsetneq \E_1\subsetneq\cdots\subsetneq \E_m=\E\subsetneq\E_{m+1}\subsetneq \cdots \subsetneq \E_{m+n}=\F.
\]
By Corollaries~\ref{87} and \ref{30}, and Remark \ref{4}.\eqref{4a}, applying $\Gamma$ yields a chain in $\GSI_S$
\[
\Gamma_{\C_\E}=\Gamma_{\E_0}\subsetneq \Gamma_{\E_1}\subsetneq\cdots\subsetneq \Gamma_{\E_m}=\Gamma_\E\subsetneq\Gamma_{\E_{m+1}}\subsetneq \cdots\subsetneq \Gamma_{\E_{m+n}}=\Gamma_\F.
\]
By Remarks~\ref{102}.\eqref{102c} and \ref{4}, Propositions~\ref{50} and \ref{62}, and Lemma~\ref{98}, dualizing with $K$ yields again a chain in $\GSI_S$
\begin{gather}\label{99}
\Gamma_{\K:\C_\E}=\Gamma_\K-\Gamma_{\C_\E}=K-\Gamma_{\E_0}\supsetneq\cdots\supsetneq K-\Gamma_{\E_m}=\\
\nonumber K-\Gamma_\E\supsetneq K-\Gamma_{\E_{m+1}}\supsetneq\cdots\supsetneq K-\Gamma_{\E_{m+n}}=K-\Gamma_\F\supset\Gamma_{\K:\F},
\end{gather}
and $d(K-\Gamma_{\E_i}\backslash K-\Gamma_{\E_{i+1}})\ge 1$ for all $i=0,\dots,m+n-1$ by Proposition~\ref{37}.
Applying Lemma~\ref{38} to the chain~\eqref{99}, it follows with equation~\eqref{140} that 
\begin{equation}\label{94}
d(K-\Gamma_{\E_i}\backslash K-\Gamma_{\E_{i+1}})=1
\end{equation}
for all $i=0,\dots,m+n-1$.
Hence,
\[
d(K-\Gamma_\E\backslash K-\Gamma_\F)=\sum_{i=m}^{m+n-1} d(K-\Gamma_{\E_i}\backslash K-\Gamma_{\E_{i+1}})=n=d(\Gamma_\F\backslash \Gamma_\E),
\]
and $d(K-\Gamma_\F\backslash \Gamma_{\K:\F})=0$.
By Proposition~\ref{37}, this implies $\Gamma_{\K:\F}=\Gamma_\K-\Gamma_\F$. 
\end{proof}


To conclude, we extend first Delgado's (see \cite[(2.8)~Thm.]{Del88}) and then Pol's (see \cite[Thm.~5.2.1]{Pol16}) characterizations of Gorensteinness to admissible rings.


\begin{dfn}
We call a good semigroup $S$ \emph{symmetric} if $S$ is a canonical ideal of $S$.
\end{dfn}


\begin{prp}\label{96}
Let $R$ be an admissible ring. 
Then $R$ is Gorenstein if and only if $\Gamma_R$ is symmetric. 
\end{prp}

\begin{proof}
Gorensteinness of $R$ is equivalent to $R$ being a canonical ideal of $R$ (see \cite[Kor.~3.4]{HK71}), and hence to $\Gamma_R$ being a canonical semigroup ideal of $\Gamma_R$ by Theorem~\ref{111}.
\end{proof}


\begin{prp}
Let $R$ be an admissible ring.
Then $R$ is Gorenstein if and only if
\begin{equation}
\label{146}
\Gamma_{R:\E} = \{\alpha \in D_{\Gamma_R} \mid \Delta^{\Gamma_R} (\tau^{\Gamma_\E}-\alpha)=\emptyset\}
\end{equation}
for every $\E \in \RFI_R$.
\end{prp}
\begin{proof}
If $R$ is Gorenstein, then $\Gamma_R$ is a canonical ideal of $\Gamma_R$ by Proposition~\ref{96}. 
Hence, Lemma~\ref{48}.\eqref{48a} and Theorems~\ref{45}.\eqref{45h} and \ref{74}.\eqref{74a} yield
\[
\Gamma_{R:\E} = \Gamma_R - \Gamma_\E = K_{\Gamma_R}^0 - \Gamma_\E = \{\alpha \in D_{\Gamma_R} \mid \Delta^{\Gamma_\E} (\tau^{\Gamma_R}-\alpha)=\emptyset\}.
\]

Conversely, if equation~\eqref{146} is satisfied for every $\E \in \RFI_R$, then, in particular,
\[
\Gamma_R = \Gamma_{R:R} = \{\alpha \in D_{\Gamma_R} \mid \Delta^{\Gamma_R} (\tau^{\Gamma_R}-\alpha)=\emptyset\}=K_{\Gamma_R}^0
\]
is a canonical ideal by Proposition~\ref{56}, and $R$ is Gorenstein by Proposition~\ref{96}.
\end{proof}

\bibliographystyle{amsalpha}
\bibliography{dsg}

\providecommand{\bysame}{\leavevmode\hbox to3em{\hrulefill}\thinspace}
\providecommand{\MR}{\relax\ifhmode\unskip\space\fi MR }
\providecommand{\MRhref}[2]{%
  \href{http://www.ams.org/mathscinet-getitem?mr=#1}{#2}
}
\providecommand{\href}[2]{#2}
\begin{thebibliography}{DdlM88}

\bibitem[AM69]{AM69}
M.~F. Atiyah and I.~G. Macdonald, \emph{Introduction to commutative algebra},
  Addison-Wesley Publishing Co., Reading, Mass.-London-Don Mills, Ont., 1969.
  \MR{0242802 (39 \#4129)}

\bibitem[BDF00]{BDF00}
V.~Barucci, M.~D'Anna, and R.~Fr\"oberg, \emph{Analytically unramified
  one-dimensional semilocal rings and their value semigroups}, J. Pure Appl.
  Algebra \textbf{147} (2000), no.~3, 215--254. \MR{1747441}

\bibitem[BH93]{BH93}
Winfried Bruns and J{\"u}rgen Herzog, \emph{Cohen-{M}acaulay rings}, Cambridge
  Studies in Advanced Mathematics, vol.~39, Cambridge University Press,
  Cambridge, 1993. \MR{MR1251956 (95h:13020)}

\bibitem[Bou61]{Bou61}
N.~Bourbaki, \emph{\'{E}l\'ements de math\'ematique. {F}ascicule {XXVII}.
  {A}lg\`ebre commutative. {C}hapitre 1: {M}odules plats. {C}hapitre 2:
  {L}ocalisation}, Actualit\'es Scientifiques et Industrielles, No. 1290,
  Herman, Paris, 1961. \MR{0217051 (36 \#146)}

\bibitem[CDK94]{CDK94}
A.~Campillo, F.~Delgado, and K.~Kiyek, \emph{Gorenstein property and symmetry
  for one-dimensional local {C}ohen-{M}acaulay rings}, Manuscripta Math.
  \textbf{83} (1994), no.~3-4, 405--423. \MR{1277539}

\bibitem[D'A97]{DAn97}
Marco D'Anna, \emph{The canonical module of a one-dimensional reduced local
  ring}, Comm. Algebra \textbf{25} (1997), no.~9, 2939--2965. \MR{1458740}

\bibitem[DdlM87]{Del87}
F\'elix Delgado de~la Mata, \emph{The semigroup of values of a curve
  singularity with several branches}, Manuscripta Math. \textbf{59} (1987),
  no.~3, 347--374. \MR{909850}

\bibitem[DdlM88]{Del88}
\bysame, \emph{Gorenstein curves and symmetry of the semigroup of values},
  Manuscripta Math. \textbf{61} (1988), no.~3, 285--296. \MR{949819}

\bibitem[Dic13]{Dic13}
Leonard~Eugene Dickson, \emph{Finiteness of the {O}dd {P}erfect and {P}rimitive
  {A}bundant {N}umbers with {$n$} {D}istinct {P}rime {F}actors}, Amer. J. Math.
  \textbf{35} (1913), no.~4, 413--422. \MR{1506194}

\bibitem[HK71]{HK71}
J{\"u}rgen Herzog and Ernst Kunz (eds.), \emph{Der kanonische {M}odul eines
  {C}ohen-{M}acaulay-{R}ings}, Lecture Notes in Mathematics, Vol. 238,
  Springer-Verlag, Berlin-New York, 1971, Seminar {\"u}ber die lokale
  Kohomologietheorie von Grothendieck, Universit{\"a}t Regensburg,
  Wintersemester 1970/1971. \MR{0412177 (54 \#304)}

\bibitem[HS06]{HS06}
Craig Huneke and Irena Swanson, \emph{Integral closure of ideals, rings, and
  modules}, London Mathematical Society Lecture Note Series, vol. 336,
  Cambridge University Press, Cambridge, 2006. \MR{2266432 (2008m:13013)}

\bibitem[J{\"a}g77]{Jag77}
Joachim J{\"a}ger, \emph{L\"angenberechnung und kanonische {I}deale in
  eindimensionalen {R}ingen}, Arch. Math. (Basel) \textbf{29} (1977), no.~5,
  504--512. \MR{0463156 (57 \#3115)}

\bibitem[Kun70]{Kun70}
Ernst Kunz, \emph{The value-semigroup of a one-dimensional {G}orenstein ring},
  Proc. Amer. Math. Soc. \textbf{25} (1970), 748--751. \MR{0265353 (42 \#263)}

\bibitem[KV04]{KV04}
K.~Kiyek and J.~L. Vicente, \emph{Resolution of curve and surface
  singularities}, Algebras and Applications, vol.~4, Kluwer Academic
  Publishers, Dordrecht, 2004, In characteristic zero. \MR{2106959
  (2005k:14028)}

\bibitem[Mat73]{Mat73}
Eben Matlis, \emph{{$1$}-dimensional {C}ohen-{M}acaulay rings}, Lecture Notes
  in Mathematics, Vol. 327, Springer-Verlag, Berlin-New York, 1973.
  \MR{0357391}

\bibitem[Mat89]{Mat89}
Hideyuki Matsumura, \emph{Commutative ring theory}, second ed., Cambridge
  Studies in Advanced Mathematics, vol.~8, Cambridge University Press,
  Cambridge, 1989, Translated from the Japanese by M. Reid. \MR{1011461
  (90i:13001)}

\bibitem[Nag62]{Nag62}
Masayoshi Nagata, \emph{Local rings}, Interscience Tracts in Pure and Applied
  Mathematics, No. 13, Interscience Publishers a division of John Wiley \&
  Sons\, New York-London, 1962. \MR{0155856 (27 \#5790)}

\bibitem[Pol16]{Pol16}
Delphine Pol, \emph{Singularit\'es libres, formes et r\'esidus logarithmiques},
  Th\`ese de doctorat, Universit\'e d'Angers, 2016.

\bibitem[Wal72]{Wal72}
Rolf Waldi, \emph{Wertehalbgruppe und {S}ingularit\"at einer ebenen algebroiden
  {K}urve}, Dissertation, Universit\"at Regensburg, 1972.

\bibitem[Wal00]{Wal00}
\bysame, \emph{On the equivalence of plane curve singularities}, Comm. Algebra
  \textbf{28} (2000), no.~9, 4389--4401. \MR{1772512}

\end{thebibliography}
\end{document}